\documentclass{article}

\usepackage{graphicx}

\usepackage[centertags]{amsmath} 
\usepackage{color} 
\usepackage{amssymb} 
\usepackage{wasysym} 
\usepackage{epsfig} 
\usepackage{rotating}
\usepackage{multirow}
\usepackage{multirow}
\usepackage{units}
\usepackage[table]{xcolor}
\usepackage{caption}
\usepackage[left=100pt,right=100pt,top=110pt, bottom=90pt]{geometry} 

\newtheorem{thm1}{Theorem}

\newtheorem{thm7}{Theorem}
\newtheorem{thm8}{Theorem}

\newtheorem{lemma}[thm1]{Lemma}

\newtheorem{fact}[thm7]{Fact}
\newtheorem{example}[thm8]{Example}

\newcommand{\R}{{\mathbb R}} 
\newcommand{\N}{{\mathbb N}} 
\newcommand{\grad}{\nabla} 
\newcommand{\lb}[1]{\underline{#1}} 
\newcommand{\ub}[1]{\overline{#1}} 
\newcommand{\hessian}[1]{{\grad^2{#1}}}

\newcommand{\func}{{\varphi}} 
\newcommand{\Order}{{\cal O}}

\newcommand{\HessianSet}{{\cal H}}
 
\newcommand{\lambdaaaT}{\Lambda_s} 
\newcommand{\lambdaabba}{\Lambda_t} 

\newcommand{\GradIntervalmulByConstShort}{c\,[y^\prime_i]}

\newcommand{\EvIntervaladdShort}{[\lambda_i]+ [\lambda_j]}

\newcommand{\EvIntervalmulByConstShort}{c\,[\lambda_i]}

\newcommand{\EvIntervalcubeShort}{3\,[y_j]\,(2\,[\lambdaaaT([y^\prime_j])]+[y_j]\,[\lambda_j])}
\newcommand{\EvIntervalpowNShort}{m[y_j]^{m-2}((m\!-\!1)[\lambdaaaT([y^\prime_j])]+[y_j]\,[\lambda_j])}

\newcommand{\GradIntervaloneOverShort}{-[y_k]^2\,[y^\prime_j]}

\newcommand{\EvIntervalexpShort}{[y_k]\,([\lambdaaaT([y^\prime_j])]+[\lambda_j])}

\newcommand{\GradIntervalsquareShort}{2\,[y_j]\,[y^\prime_j]}
\newcommand{\GradIntervalcubeShort}{3\,[y_j]^2\,[y^\prime_j]}
\newcommand{\GradIntervalpowNShort}{m\,[y_j]^{m-1}\,[y^\prime_j]}
\newcommand{\EvIntervaladdConstShort}{[\lambda_i]}

\newcommand{\GradIntervalexpShort}{[y_k]\,[y^\prime_j]}

\newcommand{\EvIntervallnShort}{1/[y_j]\,([\lambda_j]-1/[y_j]\,[\lambdaaaT([y^\prime_j])])}
\newcommand{\GradIntervalsqrtShort}{1/(2\,[y_k]) \,[y^\prime_j]}
\newcommand{\GradIntervaladdShort}{[y^\prime_i]+ [y^\prime_j]}

\newcommand{\EvIntervalmulShort}{[y_j]\,[\lambda_i]+ [y_i]\,[\lambda_j]+ [\lambdaabba([y^\prime_i], [y^\prime_j])]}

\newcommand{\EvIntervaloneOverShort}{[y_k]^2\,(2\,[y_k]\,[\lambdaaaT([y^\prime_j])]-[\lambda_j])}
\newcommand{\GradIntervallnShort}{1/[y_j]\,[y^\prime_j]}

\newcommand{\GradIntervaladdConstShort}{[y^\prime_i]}

\newcommand{\GradIntervalmulShort}{[y_i]\,[y^\prime_j]+ [y_j]\,[y^\prime_i]}

\newcommand{\EvIntervalsqrtShort}{1/(2\,[y_k])([\lambda_j]+1/(-2\,[y_j])[\lambdaaaT([y^\prime_j])])}
\newcommand{\EvIntervalsquareShort}{2\,([\lambdaaaT([y^\prime_j])]+ [y_j]\,[\lambda_j])}

\newcommand{\IAaddConstNumOps}{2\,}
\newcommand{\IAaddNumOps}{2\,}

\newcommand{\IAmulNumOps}{8\,}
\newcommand{\IAoneOverNumOps}{2\,}
\newcommand{\IAsquareNumOps}{5\,}
\newcommand{\IAcubeNumOps}{2\,}
\newcommand{\IAsqrtNumOps}{2\,}
\newcommand{\IAlnNumOps}{2\,}
\newcommand{\IAexpNumOps}{2\,}

\newcommand{\GradIntervaladdConstNumOps}{0}

\newcommand{\GradIntervaladdNumOps}{2\,n}
\newcommand{\GradIntervalmulNumOps}{18\,n}
\newcommand{\GradIntervaloneOverNumOps}{8\,n+7\,}
\newcommand{\GradIntervalsquareNumOps}{8\,n+2\,}

\newcommand{\GradIntervalsqrtNumOps}{8\,n+4\,}
\newcommand{\GradIntervalexpNumOps}{8\,n}
\newcommand{\GradIntervallnNumOps}{8\,n+2\,}

\newcommand{\EvIntervaladdConstNumOps}{0}
\newcommand{\EvIntervaladdNumOps}{2\,}
\newcommand{\EvIntervalmulNumOps}{18\,n+21\,}

\newcommand{\EvIntervalsquareNumOps}{4\,n+11\,}
\newcommand{\EvIntervalsqrtNumOps}{4\,n+25\,}

\newcommand{\EvIntervalexpNumOps}{4\,n+9\,}
\newcommand{\EvIntervallnNumOps}{4\,n+21\,}

\newcommand{\hyperRec}{B} 
\newcommand{\openDomain}{U}
\newcommand{\signMatrix}{S} 

\newcommand{\NA}{N_{\mathrm{A}}(\func)}
\newcommand{\NGersh}{N_{\mathrm{G}}(\func)}

\newcommand{\DeltaNA}{\Delta N_{\mathrm{A}}(\func)}

\newcommand{\numOfEx}{1522} 
\newcommand{\relNumPlusPlusUpper}{5.52} 
\newcommand{\relNumPlusUpper}{10.63} 
\newcommand{\relNumCircUpper}{60.33} 
\newcommand{\relNumMinusUpper}{23.52} 
\newcommand{\relNumPlusPlusLower}{2.50} 
\newcommand{\relNumPlusLower}{12.48} 
\newcommand{\relNumCircLower}{60.85} 
\newcommand{\relNumMinusLower}{24.17} 
\newcommand{\numAtLeastOneBetter}{854} 
\newcommand{\relNumAtLeastOneBetter}{56.11} 
\newcommand{\boundAllWorse}{1394} 
\newcommand{\numAllWorse}{129} 
\newcommand{\numAllBetter}{40} 
\newcommand{\meanComplexityNAOverNG}{37.40} 
\newcommand{\minComplexityNAOverNG}{10.45} 
\newcommand{\maxComplexityNAOverNG}{60.06} 
\newcommand{\meanComplexityDeltaNAOverNG}{17.35} 
\newcommand{\minComplexityDeltaNAOverNG}{3.26} 
\newcommand{\maxComplexityDeltaNAOverNG}{29.96} 
\newcommand{\rankMinComplexityNAOverNG}{18} 
\newcommand{\rankMaxComplexityNAOverNG}{161} 
\newcommand{\rankMinComplexityDeltaNAOverNG}{14} 
\newcommand{\rankMaxComplexityDeltaNAOverNG}{838} 
\newcommand{\rankFirstBreakLower}{209} 
\newcommand{\rankFirstBreakUpper}{213}

\newcommand{\approxMeanComplexityNAOverNGForNTwo}{54} 
\newcommand{\approxMeanComplexityNAOverNGForNTen}{14} 
\newcommand{\approxMeanComplexityDeltaNAOverNGForNTwo}{24} 
\newcommand{\approxMeanComplexityDeltaNAOverNGForNTen}{ 6}

\begin{document}

 {
\begin{center}
\Large
\bf
Efficient Computation of Spectral Bounds \\
for Hessian 
Matrices on Hyperrectangles \\
for Global Optimization\footnote{Funding by Deutsche Forschungsgemeinschaft is gratefully acknowledged (MO 1086/9-1).} 
\end{center}
}

 {
\begin{center}
\bf
Moritz Schulze Darup,       Martin Kastsian, Stefan Mross\\ and
Martin M\"onnigmann\footnote{M. Schulze Darup, M. Kastsian, S. Mross, M. M\"onnigmann (corr.\ author), 
              Ruhr-Universit\"at Bochum, Automatic Control and Systems Theory, 44801 Bochum, Germany.\\
        e-mail:      \texttt{martin.moennigmann@rub.de} }
\end{center}
}

\textbf{Abstract.} We compare two established and a new method for the calculation of spectral bounds for Hessian matrices on hyperrectangles by applying them to a large collection of \numOfEx{} objective and constraint functions extracted from benchmark global optimization problems. Both the tightness of the spectral bounds and the computational effort are assessed. Specifically, we compare eigenvalue bounds obtained with the interval variant of Gershgorin's circle criterion \cite{Adjiman1998a,Gershgorin1931}, Hertz and Rohn's \cite{Hertz1992,Rohn1994} method for tight bounds of interval matrices, and a recently proposed Hessian matrix eigenvalue arithmetic \cite{Monnigmann2011}, which deliberately avoids the computation of interval Hessians.\\[0.5cm] 
 

\section{Introduction}
\label{sec:Introduction}
We compare a recently proposed method \cite{Monnigmann2008a} for the calculation of spectral bounds for Hessian matrices on hyperrectangles to existing ones. 
We begin with a concise problem statement. Let 
$\func:\openDomain\subseteq\R^n\rightarrow\R$ be a twice continuously differentiable function on an open set $\openDomain\subseteq\R^n$ and let
$\hyperRec= [\lb{x}_1, \ub{x}_1]\times\cdots\times[\lb{x}_n, \ub{x}_n]$ be a closed hyperrectangle in $\openDomain$. The problem of interest reads as follows. 
\begin{equation}
\label{eq:EigenboundsOnBox}
  \begin{split}
    &\mbox{Find } \lb{\lambda}\in\R, \ub{\lambda}\in\R\mbox{ such that}\\
    &\lb{\lambda}\le \lambda\le\ub{\lambda}
    \mbox{ for all eigenvalues $\lambda$ of all matrices }
    H\in \HessianSet(\func, \hyperRec), 
  \end{split}
\end{equation}
where $\HessianSet(\func, \hyperRec)$ is the set of Hessian matrices of $\func$ on $\hyperRec$
\begin{equation}
  \label{eq:HessianSet}
  \HessianSet(\func, \hyperRec)= \left\{
      \hessian{\func}(x)\left|
        x\in \hyperRec
      \right.
    \right\}.
\end{equation}
A bound $\ub{\lambda}$ (resp.\ $\lb{\lambda}$) is called {\it tight} if there exists at least one matrix $H$ in the matrix set with
an eigenvalue $\lambda= \ub{\lambda}$ (resp.\ $\lambda= \lb{\lambda}$). Note that the bounds $\lb{\lambda}$, $\ub{\lambda}$ in (\ref{eq:EigenboundsOnBox}) may or may not be tight. 

Problem (\ref{eq:EigenboundsOnBox}) appears in various applications. It is crucial, for example, to establish the convexity of nonlinear functions in nonlinear optimization, since methods for solving nonconvex optimization problems are much less efficient than those for their convex counterparts. 
If (\ref{eq:EigenboundsOnBox}) results in $\lb{\lambda}\ge 0$ then $\func$ is convex on the
interior of the hyperrectangle $\hyperRec$ \cite{Boyd2004,Rockafellar1997}.
If, in contrast, $\lb{\lambda}<0$ results from (\ref{eq:EigenboundsOnBox}), then $\func(x)$ may or may not be convex on $\hyperRec$, but  
\begin{equation}\label{eq:alphaBB}
 \breve{\func}(x)= \func(x)- \frac{1}{2}\lb{\lambda}\sum_{i=1}^n \left(\lb{x}_i- x_i\right)\left(\ub{x}_i-x_i\right)
\end{equation}
is a convex function that underestimates $\func$ on $\hyperRec$ and coincides with $\func$ at the vertices of $\hyperRec$ \cite{Adjiman1998a,Adjiman1996}. 
Underestimators of this type are employed in nonconvex global optimization to bound the global minimum from below. 
Essentially, $\hyperRec$ is bisected into smaller and smaller hyperrectangles $\tilde{\hyperRec}$ in these approaches to obtain tighter and tighter convex underestimators. This requires solving (\ref{eq:EigenboundsOnBox}) repeatedly for different domains $\tilde{\hyperRec} \subset \hyperRec$ but the same function $\func$. 
As a result, a considerable fraction of the total computational time is spent on the calculation of convex underestimators~\cite{Adjiman1998b}. 
Consequently, fast methods for solving \eqref{eq:EigenboundsOnBox} are of interest in this field.  
Problem
(\ref{eq:EigenboundsOnBox}) also arises in automatic control
and systems theory. We refer to \cite{Monnigmann2008b} for a simple example, 
where eigenvalue bounds for Hessian matrix sets are used
to prove the positive or negative invariance of regions in the state space of
nonlinear dynamical systems.

Problem (\ref{eq:EigenboundsOnBox}) is commonly solved in two steps: 
(i) The {\it interval Hessian matrix} is calculated. 
(ii) One out of several existing methods that provide bounds on the eigenvalues of {\it symmetric interval matrices} \cite{Adjiman1998a,Hertz1992,Rohn1994} is applied. 
Interval Hessian matrices can efficiently be computed 
by combining interval
arithmetics (IA for short; see, e.g., \cite{Neumaier2008}) and automatic differentiation (see, e.g., \cite{Rall1981,Fischer1995}). 
This results in intervals $[\lb{\hessian{\func}}_{ij}, \ub{\hessian{\func}}_{ij}]\subset\R$, $i= 1, \dots, n$, $j=1, \dots, n$, such that   
\begin{equation}\label{eq:IntervalHessianElements}
  \left(\hessian{\func}(x)\right)_{ij}\in [\lb{\hessian{\func}}_{ij}, \ub{\hessian{\func}}_{ij}]  
\end{equation}
for all $x\in \hyperRec$, where $\lb{\hessian{\func}}_{ij}= \lb{\hessian{\func}}_{ji}$ and
$\ub{\hessian{\func}}_{ij}= \ub{\hessian{\func}}_{ji}$ due to symmetry of $\hessian{\func}(x)$. 
We refer to the set of matrices 
\begin{equation}\label{eq:IntervalHessian}
  \HessianSet^{\rm IA}(\func, \hyperRec) = \left\{
    H\in\R^{n\times n}\left|
      H_{ij}\in [\lb{\hessian{\func}}_{ij}, \ub{\hessian{\func}}_{ij}], \, H=H^T 
    \right.
  \right\}
\end{equation}
as the {\it interval Hessian} of $\func$ on $\hyperRec$. 
After calculating $\HessianSet^{\rm IA}(\func, \hyperRec)$, 
the spectral bounds can be found by solving the following problem.
\begin{equation}
  \label{eq:EigenboundsIntervalHessian}
  \begin{split}
    &\mbox{Find } \lb{\lambda}\in\R, 
     \ub{\lambda}\in\R \mbox{ such that}\\
    &\lb{\lambda}\le\lambda\le\ub{\lambda}
     \mbox{ for all eigenvalues $\lambda$ of all matrices }
    H\in \HessianSet^{\rm IA}(\func, \hyperRec).  
  \end{split}
\end{equation}
The calculation of $\HessianSet^{\rm IA}(\func, \hyperRec)$ requires 
$\Order(n^2)\,N(\func)$ operations if the forward mode
of automatic differentiation \cite{Fischer1995} is used, where
$N(\func)$ denotes the number of operations needed 
to evaluate $\func$ at a point in its domain. 
With the backward mode of automatic differentiation,  this complexity can be reduced
to $\Order(n)\,N(\func)$ \cite{Fischer1995}.

There exist a number of approaches to solving (\ref{eq:EigenboundsIntervalHessian}). Assuming the interval Hessian $\HessianSet^{\rm IA}(\func, \hyperRec)$ is available, the computational complexity of these methods varies between $\Order(n^2)$ for the interval variant of Gershgorin's circle criterion \cite{Adjiman1998a,Gershgorin1931} and $\Order(2^{n}\,n^3)$ for Hertz and Rohn's method \cite{Hertz1992,Rohn1994}, which provides \textit{tight} spectral bounds for  $\HessianSet^{\rm IA}(\func, \hyperRec)$ (see Sect. \ref{subSec:ClassicalApproaches} and \ref{subSec:Complexity} for details).
However, since $\HessianSet(\func, \hyperRec)\subseteq \HessianSet^{\rm IA}(\func, \hyperRec)$, problem (\ref{eq:EigenboundsIntervalHessian}) is conservative compared to the original problem (\ref{eq:EigenboundsOnBox}).
In \cite{Monnigmann2011}, we introduced a method for solving (\ref{eq:EigenboundsOnBox}) that does not require the interval Hessian $\HessianSet^{\rm IA}(\func, \hyperRec)$ and therefore avoids the conservatism inherent in  (\ref{eq:EigenboundsIntervalHessian}). The major advantage of this method is the low computational complexity, which was shown to be of order $\Order(n)\,N(\func)$ \cite{Monnigmann2011}. Note that the total numerical effort of the approaches mentioned before is the sum of the complexity for calculating $\HessianSet^{\rm IA}(\func, \hyperRec)$ \textit{and} solving (\ref{eq:EigenboundsIntervalHessian}). At least $\Order(n)\,N(\func)+\Order(n^2)$ operations are needed in these cases, and typical implementations based on forward mode automatic differentiation and Gershgorin's circle criterion require $\Order(n^2)\,N(\func)+ \Order(n^2)$ operations. 

It is the purpose of this paper to compare spectral bounds obtained with the recently proposed method \cite{Monnigmann2011} to those calculated by applying Gershgorin's circle criterion and Hertz and Rohn's method to the interval Hessian. Since the motivation for developing a new method was the application of Hessian eigenvalue bounds in global optimization, we apply the three compared methods to a large set of test functions generated from a collection of benchmark global optimization problems. Specifically, we extract \numOfEx{} objective and constraint functions 
from the COCONUT collection \cite{Shcherbina2003}. For each function, we randomly generate 100 hyperrectangles in its domain and compute the associated lower and upper eigenvalue bounds with the three methods. We compare both the resulting spectral bounds and the number of operations required by each of the approaches.

After introducing some notation in the remainder of this section, the three methods  
are summarized in Sect. \ref{sec:Arithmetic}. 
The central benchmark, which constitutes the main result of the paper, is stated in Sect.  \ref{sec:Benchmark}.
Finally, conclusions are given in Sect.  \ref{sec:Conclusion}.

\textbf{Notation.} 
Pairs of lower and upper bounds such as $\lb{\lambda}\le \lambda\le \ub{\lambda}$ are denoted by 
intervals, i.e.\ $\lambda\in[\lb{\lambda}, \ub{\lambda}]\subset\R$, for short.   
Intervals $[\lb{a}, \ub{a}]$ are further abbreviated by 
$[a]:=[\lb{a}, \ub{a}]$. Interval equality
$[a]=[b]$ is understood as $\lb{a}= \lb{b}$ and 
$\ub{a}= \ub{b}$. Calculations with intervals are carried out with
standard interval arithmetics rules, which are collected in Fact
\ref{fact:IANotation} without proof (see, e.g., \cite{Neumaier2008}). 
\begin{fact}\label{fact:IANotation}(basic interval arithmetics) 
  Let $[a]$ and $[b]$ be intervals and 
  $a\in[a]$, $b\in[b]$ and $c\in\R$ be arbitrary real
  numbers.
Then   
\begin{eqnarray} 
  a+ b &\in& [a]+ [b]:= [\lb{a}+ \lb{b}, \ub{a}+ \ub{b}],  
  \label{eq:NotationIAadd}\\ 
  a\,b &\in& [a]\,[b]:= [ 
    \min\left( 
      \lb{a}\,\lb{b}, \lb{a}\,\ub{b}, \ub{a}\,\lb{b}, \ub{a}\,\ub{b} 
    \right), 
    \max\left( 
      \lb{a}\,\lb{b}, \lb{a}\,\ub{b}, \ub{a}\,\lb{b}, \ub{a}\,\ub{b} 
    \right) 
  ]. \label{eq:NotationIAmul} \\
\label{eq:NotationIAdiv} 
  1/b&\in& 1/[b]:= [1/\ub{b},1/\lb{b}]  \\
  a+c&\in& [a]+c:= [\lb{a}+c,\ub{a}+c]  \\
 c\,a &\in& c\,[a]:= \left\{ 
    \begin{tabular}{ll}
      $[c\,\lb{a}, c\,\ub{a}]$ & if $c\ge 0$\\ 
      $[c\,\ub{a}, c\,\lb{a}]$ & if $c< 0$,
    \end{tabular}\right.
    \label{eq:NotationIAmulByConst}
\end{eqnarray} 
where $0\notin[\lb{b}, \ub{b}]$ is assumed (\ref{eq:NotationIAdiv}). 
Furthermore, the power of natural numbers $m\in\N$, the
square root, the exponential and the natural logarithm of an interval are defined as follows.  
\begin{eqnarray} 
  a^m&\in& [a^m]:= \left\{ 
    \begin{tabular}{ll} 
      $[\lb{a}^m,\ub{a}^m]$ & if $\lb{a}>0$ or $m$ odd\\
      $[\ub{a}^m,\lb{a}^m]$ & if $\ub{a}<0$ and $m$ even\\
      $[0, \max\left(-\lb{a}, \ub{a}\right)^m]$ & if $0\in [a]$ and $m$ even
    \end{tabular} 
  \right. 
  \label{eq:NotationIAsquare} \\ 
  \sqrt{a}&\in&\left[\sqrt{[a]}\right]:= [\sqrt{\lb{a}}, \sqrt{\ub{a}}],
  \label{eq:NotationIAsqrt} \\
  \exp(a)&\in& \left[\exp([a])\right]:= [\exp(\lb{a}), \exp(\ub{a})],  
  \label{eq:NotationIAexp} \\
  \ln(a)&\in&\left[\ln([a])\right]:= [\ln{\lb{a}}, \ln{\ub{a}}],
  \label{eq:NotationIAln}
\end{eqnarray} 
where $\lb{a}\ge 0$ are assumed in (\ref{eq:NotationIAsqrt}) and (\ref{eq:NotationIAln}), respectively. 
\end{fact}
By a slight abuse of notation we denote both 
a real interval $[x]= [\lb{x}, \ub{x}]\subset\R$ and a hyperrectangle $[x]= [\lb{x}, \ub{x}]=
[\lb{x}_1,  \ub{x}_1]\times\cdots\times[\lb{x}_n, \ub{x}_n]\subset\R^n$, $n\ge
2$ by a lower case letter surrounded by brackets.
As a generalization of 
Eqs.~(\ref{eq:NotationIAsquare})--(\ref{eq:NotationIAln}),   
interval extensions of functions $f(x)$, $f:\openDomain\subseteq\R^n\rightarrow\R$,
$n\ge 1$, are denoted by $[f([x])]$. 
We denote gradients and the Hessian matrices of a function $f:\openDomain\subseteq\R^n\rightarrow\R$  
by $\grad{f}(x)$ and $\hessian{f}(x)$, respectively, if they
exist. Whenever 
$\lb{\grad{f}}_i$, $\ub{\grad{f}}_i\in\R$,
$i= 1, \dots, n$ are known, then these bounds define
an interval vector denoted by $[\grad{f}]=[\lb{\grad{f}},\ub{\grad{f}}]$. Lower and upper bounds
$\lb{\hessian{f}}_{ij}$, $\ub{\hessian{f}}_{ij}\in\R$, $i= 1, \dots, n$ and
$j= 1, \dots, n$ define an interval matrix of the type (\ref{eq:IntervalHessian}), 
which is denoted by $[\hessian{f}]=[\lb{\hessian{f}},\ub{\hessian{f}}]$. Intervals
vectors and matrices are added component by component. The multiplication of an 
interval vector or matrix by an interval is understood componentwise. 
Finally, 
let $e^{(i)}\in\R^n$ be defined by $e^{(i)}_j= \delta_{ij}$, where $\delta_{ij}$ is
Kronecker's $\delta$, and let $Z$ denote the zero matrix of dimension ${n\times n}$. 

\section{Numerical calculation of eigenvalue bounds of Hessian matrices on hyperrectangles}
\label{sec:Arithmetic}

In this section we introduce the methods for the calculation of eigenvalue bounds that are applied to the collection of test cases in Sect.~\ref{sec:Benchmark}. We give only a short introduction, since these methods have been explained in detail elsewhere \cite{Adjiman1998a,Gershgorin1931,Hertz1992,Monnigmann2008a,Monnigmann2011}. 

The compared methods have in common that they are based on a codelist. A codelist results if a function $\func$ is broken down into a sequence of elementary unary and binary operations. More specifically, let $\func:\openDomain\subseteq\R^n\rightarrow\R$ denote a twice continuously differentiable function. Assume $\func$ can be evaluated at an arbitrary point $x\in\openDomain$ by carrying out a finite sequence of operations of the form 
\begin{equation}
\label{eq:CodeList} 
  \begin{array}{rcl}
        y_1 &=& x_1 \\ 
            &\vdots& \\ 
        y_n &=& x_n \\ 
    y_{n+1} &=& \Phi_{n+1} (y_1, \dots, y_n) \\ 
    y_{n+2} &=& \Phi_{n+2} (y_1, \dots, y_n, y_{n+1}) \\ 
            &\vdots& \\  
    y_{n+t} &=& \Phi_{n+t} (y_1, \dots, y_n, y_{n+1}, \dots, y_{n+t-1}) \\ 
    \func &=& y_{n+t} 
  \end{array}
\end{equation}
where each $\Phi_{n+k}$, $k=1,\dots,t$, represents one of the elementary operations listed in the first column of Tab.~\ref{tab:intervalHessianArithmetic}. The codelist lines $y_k$ for the example $\func(x_1, x_2, x_3)= \exp(x_1- 2\,x_2^2+3\,x_3^3)$ are given in the second column in (\ref{eq:CodelistExample}) below. The interval extension of a function $\func$ can be evaluated by replacing the operations in each line of (\ref{eq:CodeList}) by their interval variants listed in Fact~\ref{fact:IANotation}. For the example this results in replacing the $y_k$ from the second column by the $[y_k]$ from the third column of (\ref{eq:CodelistExample}). 

\begin{table}[h] 
\caption{Rules for the calculation of $y_k$, $[y_k]$, $[y_k^\prime]$ and $[y_k^{\prime\prime}]$ in the $k$-th line of the codelist (\ref{eq:CodeList}). $[y_k^\prime]$ refers to the interval gradient of line $k$ with respect to $x$. The operations $\Phi_k$ shown here have been selected more or less arbitrarily to accommodate a reasonably large collection of examples $\func$ treated in Sect.~\ref{sec:Benchmark}. The list of $\Phi_k$ can easily be extended \cite{Monnigmann2011}. }
\label{tab:intervalHessianArithmetic}
\centering
\begin{tabular}{lllll} 
  {\tt op} $\Phi_k$& $y_k$        & $[y_k]$                  
                                  & $[y_k^\prime]$ 
                                  & $[y_k^{\prime\prime}]$ \\\hline \hline 
  {\tt var}        & $x_k$        & $[x_k]$                      
                                  & $[e^{(k)},e^{(k)}]$
                                  & $[Z, Z]$ \\
\hline 
  {\tt addConst}   & $y_i+ c$     & $[y_i]+ [c, c]$      
                                  & $\GradIntervaladdConstShort$ 
                                  & $[y_i^{\prime\prime}]$ \\
  {\tt mulByConst} & $c\,y_i$     & $c\,[y_i]$       
                                  & $\GradIntervalmulByConstShort$
                                  & $c\,[y_i^{\prime\prime}]$ \\
  {\tt add}        & $y_i+ y_j$   & $[y_i]+ [y_j]$         
                                  & $\GradIntervaladdShort$
                                  & $[y_i^{\prime\prime}]+[y_j^{\prime\prime}]$  \\ 
  {\tt mul}        & $y_i\,y_j$   & $[y_i]\,[y_j]$           
                                  & $\GradIntervalmulShort$
                                  & $[y_i]\,[y_j^{\prime\prime}]+[y_j]\,[y_i^{\prime\prime}]+[y_i^\prime]\,[y_j^\prime]^T+[y_j^\prime]\,[y_i^\prime]^T$ \\ 
\hline 
  {\tt oneOver}    & $1/y_j$      & $1/[y_j]$      
                                  & $\GradIntervaloneOverShort$ 
                                  & $[y_k]^2\,(2\,[y_k]\,([y_j^\prime]\,[y_j^\prime]^T) -[y_j^{\prime\prime}])$ \\ 
  {\tt square}     & $y_j^2$      & $[y_j]^2$          
                                  & $\GradIntervalsquareShort$
                                  & $2\,([y_j^\prime]\,[y_j^\prime]^T+[y_j]\,[y_j^{\prime\prime}])$ \\ 
  {\tt cube}     & $y_j^3$      & $[y_j]^3$          
                                  & $\GradIntervalcubeShort$
                                  & $3\,[y_j]\,(2\,([y_j^\prime]\,[y_j^\prime])^T+[y_j]\,[y_j^{\prime\prime}])$ \\ 
  {\tt powNat}     & $y_j^m$      & $[y_j]^m$          
                                  & $\GradIntervalpowNShort$
                                  & $m\,[y_j]^{m-2}\,((m-1)\,([y_j^\prime]\,[y_j^\prime]^T)+[y_j]\,[y_j^{\prime\prime}])$ \\ 
  {\tt sqrt}       & $\sqrt{y_j}$ & $[\sqrt{[y_j]}]$           
                                  & $\GradIntervalsqrtShort$
                                  & $1/(2\,[y_k])([y_j^{\prime\prime}]+1/(-2\,[y_j])\,([y_j^\prime]\,[y_j^\prime]^T))$ \\ 
  {\tt exp}        & $\exp(y_j)$  & $[\exp([y_j])]$            
                                  & $\GradIntervalexpShort$
                                  & $[y_k]\,([y_j^\prime]\,[y_j^\prime]^T+[y_j^{\prime\prime}])$ \\ 
  {\tt ln}         & $\ln(y_j)$   & $[\ln([y_j])]$           
                                  & $\GradIntervallnShort$
                                  & $1/[y_j]\,([y_j^{\prime\prime}]-1/[y_j]\,([y_j^\prime]\,[y_j^\prime]^T))$   
\end{tabular} 
\end{table}

\subsection{Interval Hessians, Hertz and Rohn's method, and Gershgorin's circle criterion} 
\label{subSec:ClassicalApproaches}
Just as for the calculation of the interval extension of a function, a codelist can be extended to calculate gradients $\nabla \func$ and Hessians $\nabla^2 \func$ and their interval extensions $[\nabla \func]$ and $[\nabla^2 \func]$ by combining automatic differentiation (see, e.g., \cite{Fischer1995,Rall1981}) and interval arithmetics (see, e.g., \cite{Neumaier2008}). The required results are summarized in the following lemma, which summarizes results from \cite{Fischer1995}. We recall that $Z$ denotes the zero matrix of dimension $n\times n$.
\begin{lemma}\cite{Fischer1995}
\label{lem:intervalHessian}
Assume $\func$ is twice continuously differentiable on $\openDomain$ and 
can be written as a codelist. Let $\hyperRec \subset \openDomain$ be a hyperrectangle.
  Then, for all $x\in \hyperRec$, we have $\func(x)\in [\func]$,
  $\grad{\func}(x)\in[\grad{\func}]$, 
  and $\hessian{\func}(x) \in [\grad^2{\func}]$, where 
  $[\func]$, $[\grad{\func}]$, and $ [\grad^2{\func}]$ are calculated 
  by the following algorithm. 
  \begin{enumerate} 
  \item For $k= 1, \dots, n$, 
    set $[y_k]= [\lb{x}_k, \ub{x}_k]$,
    $[y_k^\prime]= [e^{(k)},e^{(k)}]$, and
    set $ [y_k^{\prime\prime}]= [Z, Z]$. 
  \item For $k= n+1, \dots, n+t$,
    calculate $[y_k]$, $[y^\prime_k]$ and $[y_k^{\prime\prime}]$
    according to columns 3$-$5 of 
    Tab.~\ref{tab:intervalHessianArithmetic}, respectively. 
  \item Set $[\func]= [y_{n+t}]$, 
    $[\grad{\func}]= [y_{n+t}^\prime]$, 
    and $[\grad^2{\func}]= [y_{n+t}^{\prime\prime}]$. 
  \end{enumerate} 
\end{lemma}

We refer to a codelist (\ref{eq:CodeList}), that has been extended by additional operations for the calculation of interval extensions or derivatives, as an {\it extended codelist} for short. The codelist (\ref{eq:CodeList}) and the extended codelist that results from Lemma \ref{lem:intervalHessian} are illustrated with an example. 
\begin{example}
\label{ex:EigBoundsA}
(interval Hessian for $\exp(x_1-2\,x_2^2+3\,x_3^3)$)
Let $\hyperRec\subset\R^3$ be an arbitrary closed hyperrectangle and consider $\func:\hyperRec\rightarrow\R$, $\func(x_1, x_2, x_3)=\exp(x_1-2\,x_2^2+3\,x_3^3)$.
Lemma \ref{lem:intervalHessian} results in the following 
expressions for $[y_k]$, $[y_k^\prime]$, and $[y_k^{\prime\prime}]$, which are first stated in a table for brevity.  
The expressions for $y_k$ stated in \eqref{eq:CodelistExample} do not result from Lemma \ref{lem:intervalHessian}, but are given for illustration of the codelist (\ref{eq:CodeList}) of $\func$. 
\begin{equation}
\label{eq:CodelistExample}
\begin{tabular}{rllllcccccc|cc} 
$k$ & $y_k$ &$[y_k]$ & $[y_k^\prime]$ & $[y_k^{\prime\prime}]$  \\ 
\hline
\hline
1 &$ x_1$ & $[x_1]$  & $([1,1],[0,0],[0,0])^T$ & $[Z,Z]$\\
2 &$ x_2$&   $[x_2]$  & $([0,0],[1,1],[0,0])^T$ & $[Z,Z]$ \\
3 &$ x_3$&   $[x_3]$  & $([0,0],[0,0],[1,1])^T$ & $[Z,Z]$ \\
\hline 
4 & $y_2^2$ &  $[y_2]^2$  & $2\,[y_2]\,[y_2^\prime]$ & $2\,([y_2^\prime]\,[y_2^\prime]^T+[y_2]\,[y_2^{\prime\prime}])$\\
5 & $y_3^3$ &  $[y_3]^3$  & $3\,[y_3]^2\,[y_3^\prime]$ & $3\,[y_3]\,(2\,([y_3^\prime]\,[y_3^\prime])^T+[y_3]\,[y_3^{\prime\prime}])$  \\
6 &$-2\,y_4$& $-2\,[y_4]$ & $-2\,[y_4^\prime]$ & $-2\,[y_4^{\prime\prime}]$\\
7 &$3\,y_5$& $3\,[y_5]$ & $3\,[y_5^\prime]$ & $3\,[y_5^{\prime\prime}]$ \\
8 &$y_1+y_6$&  $[y_1]+[y_6]$ & $[y_1^\prime]+[y_6^\prime]$ & $[y_1^{\prime\prime}]+[y_6^{\prime\prime}]$  \\
9 &$y_7+y_8$&  $[y_7]+[y_8]$ & $[y_7^\prime]+[y_8^\prime]$ & $[y_7^{\prime\prime}]+[y_8^{\prime\prime}]$\\
10 & $\exp(y_9)$ &  $[\exp([y_9])]$ & $[y_{10}]\,[y_9^\prime]$ & $[y_{10}]\,([y_9^\prime]\,[y_9^\prime]^T+[y_9^{\prime\prime}])$\\
\hline 
&$\func=y_{10}$&  $[\varphi]=[y_{10}]$ & $[\grad \varphi]=[y_{10}^\prime]$ & $[\hessian \func]=[y_{10}^{\prime\prime}]$
\end{tabular} 
\end{equation}
The codelist for $\func$ of the form (\ref{eq:CodeList}) results from rewriting the second column of (\ref{eq:CodelistExample}) as $y_1= x_1$, $y_2= x_2$, $y_3= x_3$, $y_4= y_2^2, \dots, y_{10}= \exp(y_9), \func= y_{10}$. 
The extended codelist for $[y_k^{\prime\prime}]$ can be constructed by carrying out the expressions for $[y_k]$, $[y_k^\prime]$, and $[y_k^{\prime\prime}]$ and storing the results line by line, i.e.,
\begin{equation}\label{eq:HessianCodelistExample}
 \!\!\begin{array}{lcl@{\,\,\,}lcl@{\,\,\,}lcl}
  \left[y_1\right]&=& \left[x_1\right], &
  \left[y_1^\prime\right] &=& ([1, 1], [0, 0], [0, 0])^T\!\!\!, &
  \left[y_1^{\prime\prime}\right] &=& [Z, Z], \\
  & $\vdots$ &  & 
 & $\vdots$ &  &  
& $\vdots$ &    \\
  \left[y_{10}\right] &=& \left[\exp([y_9])\right],&
  \left[y_{10}^\prime\right]&=& [y_{10}]\,[y_9^\prime], &
  \left[y_{10}^{\prime\prime}\right] &=& [y_{10}]\,([y_9^\prime]\,[y_9^\prime]^T+[y_9^{\prime\prime}]), \!\!\!\!\!\!\!\!\!\!\!
 \end{array}
\end{equation}
where the interval Hessian $[\grad^2 \func]$ reads as $[y_{10}^{\prime\prime}]$ in the codelist notation. Note that the intermediate interval function values $[y_k]$ and the derivatives $[y_k^\prime]$ are needed to calculate $[\grad^2 \func]$, while the intermediate function values $y_k$ of the original codelist for $\func$ are not. 
\end{example}

After calculating the interval Hessian $[\nabla^2 \func]$ with Lemma~\ref{lem:intervalHessian}, the relaxed problem (\ref{eq:EigenboundsIntervalHessian}) can be solved with a number of methods (see \cite{Hladik2010} for an overview).
As pointed out in Sect.~\ref{sec:Introduction}, we choose Gershgorin's circle criterion for its favorable computational complexity (see Sect.~\ref{subSec:Complexity}). In addition, we apply Hertz and Rohn's method, because it provides the {\it tight} eigenvalue bounds that solve (\ref{eq:EigenboundsIntervalHessian}). 
The interval variant of Gershgorin's circle criterion and Hertz and Rohn's method are summarized in the following two lemmas. 

\begin{lemma}[interval Gershgorin \cite{Adjiman1998a,Gershgorin1931}]
\label{lem:IntervalGershgorin}
Let $[\lb{\hessian{\func}}_{ij},\ub{\hessian{\func}}_{ij}]$, $i,j=1,\dots,n$ be intervals that define a symmetric interval matrix of the form (\ref{eq:IntervalHessian}). Then 
\begin{equation} 
  \label{eq:spectrumGershgorin}
  \lb{\lambda}= 
\min_{i\in\{1,\dots,n\}} \lb{\hessian{\func}}_{ii} - r_i\,,\qquad
  \ub{\lambda}= 
\max_{i\in\{1,\dots,n\}} \ub{\hessian{\func}}_{ii} + r_i
\end{equation}
where the Gershgorin-radii $r_i$ are defined by
$r_i=\sum_{j=1,j\neq i}^n \max(-\lb{\hessian{\func}}_{ij},\ub{\hessian{\func}}_{ij})$,
solve problem (\ref{eq:EigenboundsIntervalHessian}). 
\end{lemma}
\begin{lemma}[Hertz \cite{Hertz1992} and Rohn \cite{Rohn1994}]
\label{lem:Hertz}
Let $[\lb{\hessian{\func}}_{ij},\ub{\hessian{\func}}_{ij}]$, $i,j=1,\dots,n$ be intervals that define a symmetric interval matrix of the form (\ref{eq:IntervalHessian}). 
Define the matrices $\signMatrix^{(k)}\in\R^{n\times 2^n}$ for $k= 1, \dots, n$ recursively by
$$
  \signMatrix^{(k)} = \left(\begin{array}{cccccc}
    \multicolumn{3}{c}{\signMatrix^{(k-1)}} & \multicolumn{3}{c}{\signMatrix^{(k-1)}} 
    \\
    1 & \dots & 1 & -1 & \dots & -1 
  \end{array}\right)\,, \qquad \signMatrix^{(1)} = \left(\begin{array}{rr} 1 & -1 \end{array}\right).
$$
Define the symmetric matrices $L^{(k)}\in\R^{n\times n}$ and $U^{(k)}\in\R^{n\times n}$ 
for $k=1,\dots,2^{n-1}$ according to
$$
L^{(k)}_{ij}=\left\{\begin{tabular}{ll@{}} 
$\lb{\hessian{\func}}_{ij}$ & if $i=j$ or $\signMatrix^{(n)}_{ik} \cdot \signMatrix^{(n)}_{jk}=1$ \\
$\ub{\hessian{\func}}_{ij}$ & otherwise
\end{tabular}\right.\,,\,\,
U^{(k)}_{ij}=\left\{\begin{tabular}{ll@{}} 
$\ub{\hessian{\func}}_{ij}$ & if $L^{(k)}_{ij}=\lb{\hessian{\func}}_{ij}$ \\
$\lb{\hessian{\func}}_{ij}$ & otherwise
\end{tabular}\right..
$$
Then
\begin{equation} 
  \label{eq:spectrumHertz} 
  \lb{\lambda}= 
\min_{k\in\{1,\dots,2^{n-1}\}} \lambda_{\min}(L^{(k)}) \,,\qquad
  \ub{\lambda}= 
\max_{k\in\{1,\dots,2^{n-1}\}} \lambda_{\max}(U^{(k)}),
\end{equation}
where $\lambda_{\min}(A)$ and $\lambda_{\max}(A)$ denote the smallest and largest (real) eigenvalue of any symmetric real matrix $A=A^T$, respectively, solve problem (\ref{eq:EigenboundsIntervalHessian}). 
\end{lemma}
We briefly illustrate Lemmas \ref{lem:IntervalGershgorin} and \ref{lem:Hertz} by applying them to the sample function from Example \ref{ex:EigBoundsA}.
\begin{example}
\label{ex:EigBoundsGH}
(Gershgorin and Hertz/Rohn applied to $\exp(x_1-2\,x_2^2+3\,x_3^3)$) 
Without detailing the calculations we claim that substituting   $\hyperRec=[-0.3,0.2]\times[-0.1,0.6]\times[-0.4,0.5]$ into the extended codelist \eqref{eq:HessianCodelistExample} yields
$$[\hessian \func]=[y^{\prime\prime}_{10}]=\left(\begin{array}{ccc} 
      \left[0.298,1.777\right] &[-4.265,0.7109] &   [0.000 , 3.999]       \\
      \left[-4.265,0.7109\right] &   [-7.109 , 3.128] &   [-9.597,1.599]   \\
            \left[0.000 , 3.999\right] &  [-9.597,1.599] &[ -12.795,24.991]
\end{array}\right).
$$
Applying Lemma \ref{lem:IntervalGershgorin} yields the $n=3$ Gershgorin radii
$r_1=4.265+3.999=8.264$, $r_2=4.265+9.597=13.862$ and $r_3=3.999+9.597=13.596$. Upon substitution into (\ref{eq:spectrumGershgorin}) the spectral bounds 
\begin{equation}
\label{eq:resultExampleGershgorin}
[\lambda_{\mathrm{G}}]
=[\lb{\lambda}_{\mathrm{G}},\ub{\lambda}_{\mathrm{G}}]
=[-12.795-13.596,24.991+13.596]=[-26.391,38.587]
\end{equation}
result, where the subscript G is short for Gershgorin. Hertz and Rohn's method requires to calculate $2\cdot 2^{n-1}=2^n=8$ vertex matrices $L^{(1)},\dots,L^{(4)}$ and $U^{(1)},\dots,U^{(4)}$ and the sign matrix $\signMatrix^{(3)}$ defined in Lemma \ref{lem:Hertz}.
Equation (\ref{eq:spectrumHertz}) yields 
\begin{equation}
\label{eq:resultExampleHertz}
[\lambda_{\mathrm{H}}]=[\lb{\lambda}_{\mathrm{H}},\ub{\lambda}_{\mathrm{H}}]= 
[\lambda_{\min}(L^{(2)}),\lambda_{\max}(U^{(3)})] = [-20.597, 29.603] 
\end{equation}
with $L^{(2)}=  \left(\begin{array}{rrr}  0.298  &     0.711  &      3.999 \\
     0.711   &    -7.109   &    -9.597 \\
           3.999    &   -9.597   &   -12.795 \end{array}\right)$ and $U^{(3)}=   \left(\begin{array}{rrr}  1.777     &  -4.265    &    3.999 \\
       -4.265   &     3.128    &   -9.597 \\
        3.999   &    -9.597    &   24.991 \end{array}\right)$,
where the subscript H is short for Hertz and Rohn. 
We only list the matrices $L^{(2)}$ and $U^{(3)}$ that are selected in the minimization and maximization in (\ref{eq:spectrumHertz}) and omit the remaining six vertex matrices for brevity. 
\end{example}

\subsection{Eigenvalue arithmetic} 

We summarize the eigenvalue arithmetic in Lemma \ref{lem:directArithmetic} and Tab.~\ref{tab:directArithmetic} and refer the reader to \cite{Monnigmann2008a,Monnigmann2011} for details.
Lemma~\ref{lem:directArithmetic} implies that the eigenvalue arithmetic does {\it not} require the interval Hessian, but interval gradients suffice. 
The functions $\lambdaaaT$ and $\lambdaabba$ used in Lemma~\ref{lem:directArithmetic} and Tab.~\ref{tab:directArithmetic} are defined as
\begin{equation} 
  \label{eq:lambdaaaT} 
[\lambdaaaT([a])]= \left[0,\sum_{i= 1}^n \max(\lb{a}_i^2, \ub{a}_i^2)\right]  
\end{equation}
 and
\begin{equation}\label{eq:Lambda_abbaNotation} 
  [\lambdaabba([a], [b])]= 
      [-\beta, \beta]+  
      \sum\limits_{i=1}^n [\lb{a}_i, \ub{a}_i]\,[\lb{b}_i, \ub{b}_i]
\end{equation} 
where $\beta= \sqrt{
    (\sum_{i=1}^n \max(\lb{a}_i^2, \ub{a}_i^2)) 
    (\sum_{i=1}^n \max(\lb{b}_i^2, \ub{b}_i^2)) 
  }$. We refer to  \cite{Monnigmann2011}  for a detailed discussion of the meaning of $[\lambdaaaT([a])]$ and $[\lambdaabba([a],[b])]$.

\begin{table}[h] 
\caption{Rules for the calculation of $[\lambda_k]$ in the $k$-th line of the codelist (\ref{eq:CodeList}). We assume that 
$[y_i]$, $[y_i^\prime]$ and 
$[y_j]$, $[y_j^\prime]$  
for all previous lines $i\leq k$, $j\leq k$ have been calculated according to the rules from Tab.~\ref{tab:intervalHessianArithmetic} and can be reused in line $k$. Rules for $y_k$ are repeated here for convenience. The functions $[\lambdaaaT]([a])$ and $[\lambdaabba]([a],[b])$ are defined in (\ref{eq:lambdaaaT}) and (\ref{eq:Lambda_abbaNotation}).}
\label{tab:directArithmetic}
\centering
\begin{tabular}{lllll} 
  {\tt op} $\Phi_k$& $y_k$        &  $[\lambda_k]$  \\ \hline \hline 
  {\tt var}        & $x_k$        &  $[0,0]$  \\
\hline 
  {\tt addConst}   & $y_i+ c$     &  $\EvIntervaladdConstShort$  \\
  {\tt mulByConst} & $c\,y_i$     & $\EvIntervalmulByConstShort$ \\
  {\tt add}        & $y_i+ y_j$   & $\EvIntervaladdShort$   \\ 
  {\tt mul}        & $y_i\,y_j$   & $\EvIntervalmulShort$  \\ 
\hline 
  {\tt oneOver}    & $1/y_j$      &  $\EvIntervaloneOverShort$ \\ 
  {\tt square}     & $y_j^2$      &  $\EvIntervalsquareShort$\\ 
  {\tt cube}     & $y_j^3$      &  $\EvIntervalcubeShort$ \\ 
  {\tt powNat}     & $y_j^m$      &  $\EvIntervalpowNShort$ \\ 
  {\tt sqrt}       & $\sqrt{y_j}$ &   $\EvIntervalsqrtShort$ \\ 
  {\tt exp}        & $\exp(y_j)$  &  $\EvIntervalexpShort$        \\
  {\tt ln}         & $\ln(y_j)$   &     $\EvIntervallnShort$      
\end{tabular} 
\end{table} 

\begin{lemma}[eigenvalue arithmetic \cite{Monnigmann2011}]
\label{lem:directArithmetic}
Assume $\func$ is twice continuously differentiable on $\openDomain$ and can be written as a codelist. Let $\hyperRec \subset \openDomain$ be a hyperrectangle.
  Then, for all $x\in \hyperRec$, we have $\func(x)\in [\func]$,
  $\grad{\func}(x)\in[\grad{\func}]$, 
  and $\lambda_\func\in [\lambda_\func]$ for all eigenvalues of
  $\hessian{\func}(x)$, where 
  $[\func]$, $[\grad{\func}]$, and $[\lambda_\func]$ are calculated 
  by the following algorithm.  
  \begin{enumerate} 
  \item For $k= 1, \dots, n$, 
    set $[y_k]= [\lb{x}_k, \ub{x}_k]$,
    $[y_k^\prime]= [e^{(k)},e^{(k)}]$, and
    set $[\lambda_k]= [0, 0]$. 
  \item For $k= n+1, \dots, n+t$,
    calculate $[y_k]$, $[y^\prime_k]$ and $[\lambda_k]$
    according to columns 3 and 4 of
    Tab.~\ref{tab:intervalHessianArithmetic} and column 3 of Tab.~\ref{tab:directArithmetic}, respectively.   
  \item Set $[\func]= [y_{n+t}]$, 
    $[\grad{\func}]= [y_{n+t}^\prime]$, 
    and $[\lambda_\func]= [\lambda_{n+t}]$. 
  \end{enumerate} 
\end{lemma}

Lemma \ref{lem:directArithmetic} is illustrated with 
the sample function from Example \ref{ex:EigBoundsA} and \ref{ex:EigBoundsGH}. 
\begin{example}\label{ex:NewEigboundExample}
Let $\hyperRec$ and $\func:\hyperRec\rightarrow\R$ be as in Example \ref{ex:EigBoundsA}. Applying Lemma \ref{lem:directArithmetic} to $\func$ results in the expressions for $[y_k]$, $[y_k^\prime]$, and $[\lambda_k]$ listed in \eqref{eq:NewEigboundExample}, which we state in a table first for brevity. Note that $[y_k]$ and $[y_k^\prime]$ are equal to those in \eqref{eq:CodelistExample}. These expressions are repeated here, since the $[\lambda_k]$ depend on them. 
\begin{equation}
\label{eq:NewEigboundExample}
\begin{tabular}{rlllcccccc|cc} 
$k$ &$[y_k]$ & $[y_k^\prime]$ & $[\lambda_k]$  \\ 
\hline
\hline
1 & $[x_1]$  & $([1,1],[0,0],[0,0])^T$ & $[0,0]$\\
2 &   $[x_2]$  & $([0,0],[1,1],[0,0])^T$ & $[0,0]$ \\
3 &   $[x_3]$  & $([0,0],[0,0],[1,1])^T$ & $[0,0]$ \\
\hline 
4 &  $[y_2]^2$  & $2\,[y_2]\,[y_2^\prime]$ & $2\,([\lambdaaaT([y^\prime_2])]+[y_2]\,[\lambda_2])$\\
5 &  $[y_3]^3$  & $3\,[y_3]^2\,[y_3^\prime]$ & $3\,[y_3]\,(2\,[\lambdaaaT([y^\prime_3])]+[y_3]\,[\lambda_3])$  \\
6 & $-2\,[y_4]$ & $-2\,[y_4^\prime]$ & $-2\,[\lambda_4]$\\
7 & $3\,[y_5]$ & $3\,[y_5^\prime]$ & $3\,[\lambda_5]$ \\
8 &  $[y_1]+[y_6]$ & $[y_1^\prime]+[y_6^\prime]$ & $[\lambda_1]+[\lambda_6]$  \\
9 &  $[y_7]+[y_8]$ & $[y_7^\prime]+[y_8^\prime]$ & $[\lambda_7]+[\lambda_8]$\\
10 &  $[\exp([y_9])]$ & $[y_{10}]\,[y_9^\prime]$ & $[y_{10}]\,([\lambdaaaT([y^\prime_9])]+[\lambda_9])$\\
\hline 
&  $[\varphi]=[y_{10}]$ & $[\grad \varphi]=[y_{10}^\prime]$ & $[\lambda_{\varphi}]=[\lambda_{10}]$
\end{tabular} 
\end{equation}

The extended codelist for $[\lambda_\func]$ results from evaluating and storing the expressions listed in \eqref{eq:NewEigboundExample} line by line, i.e.,
\begin{equation}\label{eq:CodelistNewEigboundExample}
 \!\begin{array}{lcl@{\,\,\,}lcl@{\,\,\,}lcl}
  \left[y_1\right]&=& \left[x_1\right], &
  \left[y_1^\prime\right] &=& ([1, 1], [0, 0], [0, 0])^T\!\!\!, &
  \left[\lambda_1\right] &=& [0, 0], \\
  & $\vdots$ &  & 
 & $\vdots$ &  &  
& $\vdots$ &    \\
  \left[y_{10}\right] &=& \left[\exp([y_9])\right],&
  \left[y_{10}^\prime\right]&=& [y_{10}]\,[y_9^\prime], &
  \left[\lambda_{10}\right] &=& [y_{10}]\,([\lambdaaaT([y^\prime_9])]+[\lambda_9]), \!\!\!\!\!\!\!\!\!\!\!
 \end{array}
\end{equation}
Without detailing the calculations we claim that applying \eqref{eq:CodelistNewEigboundExample} to the particular hyperrectangle 
$\hyperRec=[-0.3,0.2]\times[-0.1,0.6]\times[-0.4,0.5]$ from Example \ref{ex:EigBoundsGH} results in 
\begin{equation}
\label{eq:resultExampleArithmetic}
[\lambda_{\mathrm{A}}]=[\lambda_{\func}]=[-19.904,37.004],
\end{equation}
where the subscript A is short for arithmetic. 
\end{example}

By comparing the spectral bounds (\ref{eq:resultExampleGershgorin}), (\ref{eq:resultExampleHertz}) and (\ref{eq:resultExampleArithmetic}), we find the relations $\lb{\lambda}_{\mathrm{G}}<\lb{\lambda}_{\mathrm{H}}<\lb{\lambda}_{\mathrm{A}}$ and $\ub{\lambda}_{\mathrm{H}}<\ub{\lambda}_{\mathrm{A}}<\ub{\lambda}_{\mathrm{G}}$
for the discussed example. Note that the lower bound from the eigenvalue arithmetic is tighter than the tight bound for the interval Hessian obtained with Hertz and Rohn's method. 
We stress that these relations do not hold in general. It is the very point of Section \ref{sec:Benchmark} to analyze these relations for a large collection of examples.

\subsection{Computational complexities}
\label{subSec:Complexity}
The discussed methods do not only differ with respect to the tightness of the eigenvalue bounds, but also with respect to computational cost. Calculating the interval Hessian matrix with forward mode automatic differentiation and applying 
Hertz and Rohn's method
requires 
\begin{equation}
  \Order(n^2)\,N(\func)+\Order(2^n n^3)\label{eq:ComplexityHertz}
\end{equation}
operations \cite{Monnigmann2011},  
where $N(\varphi)$ denotes the number of operations needed for the evaluation of $\varphi$ at a point. Calculating the interval Hessian and applying Gershgorin's circle criterion
 takes
\begin{equation}
  \Order(n^2)\,N(\func)+\Order(n^2)\label{eq:ComplexityGershgorin}
\end{equation}
operations~\cite{Monnigmann2011}. 
Note that $\Order(n^2)\,N(\func)$ operations are needed for the calculation of the interval Hessian in both \eqref{eq:ComplexityHertz} and \eqref{eq:ComplexityGershgorin}. Calculating eigenvalue bounds with the arithmetic from \cite{Monnigmann2011} requires
\begin{equation}\label{eq:ComplexityArithmetic}
  \Order(n)\,N(\func)
\end{equation}
operations \cite{Monnigmann2011}. Due to the $\Order(2^n\,n^3)$ term in (\ref{eq:ComplexityHertz}) the computational cost of Hertz and Rohn's grows drastically compared to (\ref{eq:ComplexityGershgorin}) and (\ref{eq:ComplexityArithmetic}). The complexities (\ref{eq:ComplexityGershgorin}) and (\ref{eq:ComplexityArithmetic}), however, differ only by one order of magnitude. We therefore compare the computational effort of these two methods more precisely in the present and subsequent section. 

We denote the eigenvalue bounds calculated with the arithmetic proposed here and Gershgorin's circle criterion $[\lambda_\mathrm{A}]$ and $[\lambda_\mathrm{G}]$, respectively. The number of operations required to calculate $[\lambda_\mathrm{A}]$ and $[\lambda_\mathrm{G}]$ for a specific function $\func$ are denoted $\NA$ and $\NGersh$, respectively.

\begin{table}[h]
\caption{Number of operations necessary to calculate $y_k$, $[y_k]$, $[y^\prime_k]$, $[\lambda_k]$ and $[y^{\prime\prime}_k]$ for each type of line of a codelist (\ref{eq:CodeList}). $N([y^\prime_k])$ denotes the number of operations necessary to compute $[y^\prime_k]$ assuming that $[y_k]$ is already available. $N([\lambda_k])$ and $N([y^{\prime\prime}_k])$ denote the number of operations necessary to compute $[\lambda_k]$ and $[y^{\prime\prime}]$, respectively, assuming $[y_k]$ and $[y^\prime_k]$ are available. Listed numbers apply for $n> 1$. 
}

\label{tab:complexityOperations}
\centering
\begin{tabular}{lccccc} 
  {\tt op} $\Phi_k$& $N(y_k)$     & $N([y_k])$                   
                                  & $N([y_k^\prime])$ 
                                  & $N([\lambda_k])$ 
                                  & $N([y_k^{\prime\prime}])$ \\\hline \hline 
  {\tt var}        & $1$        & $0$                      
                                  & $0$
                                  & $0$ 
                                  & $0$  \\
\hline 
  {\tt addConst}   & $1$     & $\IAaddConstNumOps$      
                                  & $\GradIntervaladdConstNumOps$ 
                                  & $\EvIntervaladdConstNumOps$
                                  & $0$ \\
  {\tt mulByConst} & $1$     & \!$2$ &  $2\,n$ & $2$ & $n\,(n+1)$\\   
  {\tt add}        & $1$   & $\IAaddNumOps$        
                                  & $\GradIntervaladdNumOps$
                                  & $\EvIntervaladdNumOps$ 
                                  &   $n\,(n+1)$ \\ 
  {\tt mul}        & $1$   & $\IAmulNumOps$          
                                  & $\GradIntervalmulNumOps$
                                  & $\EvIntervalmulNumOps$ 
                                  &  $19\,n\,(n+1)$ \\
\hline 
  {\tt oneOver}    & $1$      & $\IAoneOverNumOps$       
                                  & $\GradIntervaloneOverNumOps$ 
                                  & $4\,n+26$
                                  & $14\,n\,(n+1) + 7$ \\
  {\tt square}     & $1$      & $\IAsquareNumOps$          
                                  & $\GradIntervalsquareNumOps$
                                  & $\EvIntervalsquareNumOps$
                                  & $10\,n\,(n+1)$ \\
  {\tt cube}     & $1$      & $\IAcubeNumOps$         
                                  & $8\,n+7$
                                  & $4\,n+21$ 
                                  & $14\,n\,(n+1) + 2$ \\
  {\tt powNat}    & 1      & 5 & $8\,n+7$ & $4\,n+26$ & $14\,n\,(n+1) + 7$  \\
  {\tt sqrt}       & $1$ & $\IAsqrtNumOps$           
                                  & $\GradIntervalsqrtNumOps$
                                  & $\EvIntervalsqrtNumOps$
                                  & $13\,n\,(n+1)+8$  \\
  {\tt exp}        & $1$  & $\IAexpNumOps$             
                                  & $\GradIntervalexpNumOps$
                                  & $\EvIntervalexpNumOps$ 
                                  & $9\,n\,(n+1)$  \\
  {\tt ln}         & $1$   & $\IAlnNumOps$            
                                  & $\GradIntervallnNumOps$
                                  & $\EvIntervallnNumOps$ 
                                  &  $13\,n\,(n+1) +4$   
\end{tabular} 
\end{table} 
The exact number of operations needed to calculate eigenvalue bounds for a specific function can be determined for any of the discussed methods by counting operations in the extended codelist of $\func$. 
Table \ref{tab:complexityOperations} lists the exact number of operations needed in each codelist line by line type. 
An operation counted towards $N(y_k)$, $N([y_k])$, $N([y_k^\prime])$, $N([\lambda_k])$, or $N([y_k^{\prime\prime}])$ in Tab.~\ref{tab:complexityOperations} may either be an addition, multiplication or comparison of two real numbers, or the application of one of the functions \texttt{oneOver}, \texttt{square}, \texttt{cube}, \texttt{powNat}, \texttt{sqrt}, \texttt{exp} or \texttt{ln}.
Note that this way of counting operations is coarse but a standard approach in the field of automatic differentiation \cite{Fischer1995,Rall1981}.

Before applying Tab.~\ref{tab:complexityOperations} to specific examples in Sect.~\ref{sec:Benchmark}, we derive some general statements. 
From the last two columns of Tab.~\ref{tab:complexityOperations} we infer $N([\lambda_k]) \leq N([y_k^{\prime\prime}])$ for all codelist line types (assuming $n>1$). Since the $[y_k^{\prime\prime}]$ are required for the calculation of eigenvalue bounds with Gershgorin's circle criterion, $N([\lambda_k]) \leq N([y_k^{\prime\prime}])$ for all $\Phi_k$ implies 
\begin{equation}
\label{eq:complexityRelationAG}
\NA
 \leq \NGersh
\end{equation} 
for any function $\func$ that can be stated as a codelist with lines of the types from Tab.~\ref{tab:complexityOperations}. 
Furthermore, inspection of Tab.~\ref{tab:complexityOperations} shows that  
the eigenvalue arithmetic can be applied at little additional computational effort, whenever eigenvalue bounds are calculated by applying Gershgorin's method to the interval Hessian. 
This statement holds, since the $[y_k^\prime]$ required for the arithmetic are available as an intermediate result to the interval Hessian calculation. 
More specifically,  
\begin{equation}\label{eq:additionalComplexity}
  \DeltaNA=  \sum \limits_{k=1}^t N([\lambda_{n+k}])
\end{equation} 
additional operations are needed to calculate eigenvalue bounds with the arithmetic, if they are calculated by applying Gershgorin's circle criterion to the interval Hessian matrix. We infer from Tab.~\ref{tab:complexityOperations} that $\Delta N_\mathrm{A}(\func)$ as defined in (\ref{eq:additionalComplexity}) amount to $\Order(n)$ operations. 

\begin{example}(number of operations for Examples~\ref{ex:EigBoundsA}, \ref{ex:EigBoundsGH} and \ref{ex:NewEigboundExample})
Table \ref{tab:ComplexityExample} lists the particular numbers of operations necessary to evaluate $[y_k]$, $[y_k^\prime]$, $[\lambda_k]$ and $[y_k^{\prime\prime}]$ for the sample function $\func(x)=\exp(x_1-2\,x_2^2+3\,x_3^3)$ according to the codelist  (\ref{eq:CodelistExample}). We find $\NA=N([\lambda_{\mathrm{A}}])=N([\lambda_\func])=17+73+73=163$ and $\NGersh=N([\hessian \func]) + N([\lambda_{\mathrm{G}}])=17+73+222+12=324$, where $N([\lambda_{\mathrm{G}}])=12$. The additional effort for calculating $\lambda_{\mathrm{A}}$ given the intermediate results $[y_k]$ and $[y_k^\prime]$ yields $\DeltaNA=73$.
\begin{table}[h]
\caption{Numbers of operations for the codelist lines of Example \ref{ex:EigBoundsA}. }
\label{tab:ComplexityExample}
\centering
\begin{tabular}{rlrrrr} 
$k$ &  {\tt op} $\Phi_k$& $N([y_k])$ & $N([y_k^\prime])$ & $N([\lambda_k])$ & $N([y_k^{\prime\prime}])$  \\ 
\hline
\hline
4 & {\tt square} & 5 & 18 & 19 & 60  \\
5 & {\tt cube} & 2 & 23 & 29 & 86 \\
6 & {\tt mulByConst} & 2 & 4 & 2 & 6 \\
7 & {\tt mulByConst} & 2 & 4 & 2 & 6   \\
8 & {\tt add} &  2 & 4 & 2 & 6  \\
9 & {\tt add} &  2 & 4 & 2 & 6   \\
10 & {\tt exp} & 2 & 16 & 17 & 54 \\
\hline 
$\sum$ & & 17 & 73 & 73 & 222
\end{tabular} 
\end{table}
\end{example}

\section{Benchmark: Arithmetic versus Gershgorin and Hertz}
\label{sec:Benchmark}
We apply the eigenvalue arithmetic to a large collection of examples and compare results to those obtained by applying Gershgorin's circle criterion and Hertz and Rohn's method to the interval Hessian. Sections \ref{subsec:Coconut} and \ref{subsec:RatingScheme} describe the test examples and the scheme of comparison. The actual results are summarized in Sect.~\ref{subsec:Results}

\subsection{Collection of test cases}
\label{subsec:Coconut}
The test cases are extracted from the COCONUT collection of optimization problems \cite{Shcherbina2003}. 
We consider all COCONUT problems with $1<n\le 10$ variables and extract those cost and constraint functions 
that can be decomposed into the operations listed in Tabs. \ref{tab:intervalHessianArithmetic} and \ref{tab:directArithmetic} respectively. This results in a set of \numOfEx{} sample functions $\func:\R^n\rightarrow\R$ with $1<n\le 10$.  
For each $\func$, we generate 100 random hyperrectangles $\hyperRec\subseteq D\subset\R^n$ in the domain $D$ of $\func$ specified in the respective COCONUT problem. Each of the three methods introduced in Lemmas~\ref{lem:IntervalGershgorin}--\ref{lem:directArithmetic} is applied to the resulting $\numOfEx\cdot 100$ sample problems.

We omit examples with $n=1$, since the three methods yield identical spectral bounds and require the same numerical effort. The upper bound $n\le 10$ is arbitrary. The comparison in Sect.~\ref{subsec:Results} corroborates that the eigenvalue arithmetic benefits more and more from its 
favorable computational complexity as $n$ increases, which was anticipated in the comparison of computational complexities in Sect.~\ref{subSec:Complexity}. While the eigenvalue arithmetic and Gershgorin's circle criterion could be applied well beyond $n= 10$, it becomes tedious to calculate the exact Hessian matrix eigenvalue bounds for comparison, due to the $\Order(2^n\,n^3)$ complexity of this problem.

\begin{table}[h]
\caption{Excerpt of the set of examples taken from the COCONUT-benchmark.}
\label{tab:examples}
\centering
\begin{tabular}{cc|c}
name & $n$ & function $\func$\\
\hline
\hline
\vspace{-2mm}
&  \\
\verb"ex8_1_6-1" & 2 & $\frac{1}{(x_1-4)^2+(x_2-4)^2+0.1}+\frac{1}{(x_1-1)^2+(x_2-1)^2+0.2}+\frac{1}{(x_1-8)^2+(x_2-8)^2+0.2}$ \\\vspace{-2mm}
&  \\
\verb"ex7_2_6-2" &3 & $1-0.01\,\frac{x_2}{x_3}-0.01\,x_1-0.0005\,x_1\,x_3$ \\
\vspace{-2mm}
&  \\
\verb"ex14_2_2-6" & 4 & $10.208-\frac{2755.642}{x_3+219.161}-\frac{0.192\,x_1}{x_1+0.192\,x_2}-\frac{x_2}{0.316\,x_1+x_2}-\ln(0.316\,x_1+x_2)+x_4$\\
\end{tabular}
\end{table}

Table \ref{tab:examples} lists three sample functions from the COCONUT collection for illustration. 
We refer to all examples by their COCONUT name, for example \verb"ex8_1_6". 
The suffix \verb"-i", as in \verb"ex8_1_6-1" for example, uniquely identifies the function in the respective COCONUT optimization problem, where $i=1$ corresponds to the objective function and $i=2,\dots,m$ corresponds to the $(i-1)$-th constraint function\footnote{Ordering is as in the GAMS code provided in the COCONUT library. }.

We stress that we use the described set of test functions without further modifications. There exist functions in the collection treated here that contain convex terms, or terms for which tight convex under- or tight concave overestimators are known (e.g., bilinear, trilinear, linear fractional terms) \cite{Adjiman1998a,McCormick1976} . 
Depending on the application it may be advisable to separate these terms from the given function $\func$, and to calculate Hessian eigenvalue bounds only for the remaining terms of $\func$. Here we choose not to apply any preprocessing for the sake of an unbiased comparison. 

\subsection{Evaluation of results}
\label{subsec:RatingScheme}

We introduce a simple rating scheme that assigns each result to one of a finite set of classes. Specifically, we distinguish between the cases listed in Tab.~\ref{tab:ratingScheme}, which reflect that a bound from the eigenvalue arithmetic may be 
\begin{equation*}
 \begin{array}{cl}
   (-) & \mbox{ worse than the bounds from the other two methods},
   \\
   (\circ) & \mbox{ equal to the one from Gershgorin's method, but equal to or worse than}\\ 
           & \mbox{ the one from Hertz and Rohn's method}, 
   \\
   (+) & \mbox{ better than the one from Gershgorin's method, and equal to or worse}\\
       & \mbox{ than the one from Hertz and Rohn's method}, 
   \\
   (++) & \mbox{ better than the one from Hertz and Rohn's method}. 
 \end{array}
\end{equation*}
Note that bounds calculated with Hertz and Rohn's method are never worse than those from Gershgorin's circle criterion, since Hertz and Rohn's method provides the tight eigenvalue bounds for an interval matrix. Consequently, the bounds from Gershgorin's method do not play a role in our definition of the $(++)$ category. Furthermore note that we do not distinguish between 
$\lb{\lambda}_{\mathrm{G}} = \lb{\lambda}_{\mathrm{H}}$ 
and $\lb{\lambda}_{\mathrm{G}} < \lb{\lambda}_{\mathrm{H}}$
(resp.\ $\ub{\lambda}_{\mathrm{G}} = \ub{\lambda}_{\mathrm{H}}$ 
and $\ub{\lambda}_{\mathrm{G}} > \ub{\lambda}_{\mathrm{H}}$)
in the case 
$\lb{\lambda}_{\mathrm{A}}\leq\lb{\lambda}_{\mathrm{G}}$
(resp.\ $\ub{\lambda}_{\mathrm{A}}\geq\ub{\lambda}_{\mathrm{G}}$). 

\begin{table}[h]
\caption{Classes used to aggregate results in Sect.~\ref{subsec:Results}. Symbols  $[\lambda_{\mathrm{A}}]=[\lb{\lambda}_{\mathrm{A}},\ub{\lambda}_{\mathrm{A}}]$,  $[\lambda_{\mathrm{G}}]=[\lb{\lambda}_{\mathrm{G}},\ub{\lambda}_{\mathrm{G}}]$, and $[\lambda_{\mathrm{H}}]=[\lb{\lambda}_{\mathrm{H}},\ub{\lambda}_{\mathrm{H}}]$ denote the eigenvalue bounds calculated with the eigenvalue arithmetic, Gershgorin's circle criterion, and Hertz and Rohn's method, respectively.}
\label{tab:ratingScheme}
\centering
\begin{tabular}{c|c|c|c|c}
 bound & \multicolumn{4}{c}{class} \\
  & ($-$) & ($\circ$) & ($+$) & ($++$)\\ 
\hline 
\hline
upper ($\ub{\lambda}_{\mathrm{A}}$) &   $\ub{\lambda}_{\mathrm{A}}>\ub{\lambda}_{\mathrm{G}}\geq\ub{\lambda}_{\mathrm{H}}$ & $\ub{\lambda}_{\mathrm{A}}=\ub{\lambda}_{\mathrm{G}}\geq\ub{\lambda}_{\mathrm{H}}$ & $\ub{\lambda}_{\mathrm{G}}>\ub{\lambda}_{\mathrm{A}}\geq\ub{\lambda}_{\mathrm{H}}$ & $\ub{\lambda}_{\mathrm{G}}\geq\ub{\lambda}_{\mathrm{H}}>\ub{\lambda}_{\mathrm{A}}$\\
lower ($\lb{\lambda}_{\mathrm{A}}$)& 
$\lb{\lambda}_{\mathrm{A}}<\lb{\lambda}_{\mathrm{G}}\leq\lb{\lambda}_{\mathrm{H}}$ & $\lb{\lambda}_{\mathrm{A}}=\lb{\lambda}_{\mathrm{G}}\leq\lb{\lambda}_{\mathrm{H}}$ & $\lb{\lambda}_{\mathrm{G}}<\lb{\lambda}_{\mathrm{A}}\leq\lb{\lambda}_{\mathrm{H}}$ & $\lb{\lambda}_{\mathrm{G}}\leq\lb{\lambda}_{\mathrm{H}}<\lb{\lambda}_{\mathrm{A}}$ 
\end{tabular}
\end{table}

Table \ref{tab:casesExamples} lists some numerical results. These examples illustrate that the classes introduced in Tab.~\ref{tab:ratingScheme} are meaningful. In particular it is evident that eigenvalue bounds of the same function may fall into different classes for different hyperrectangles $\hyperRec$.

\begin{table}[h]
\caption{Illustration of the classes $(-)$, $(\circ)$, $(+)$, $(++)$ introduced in Tab.~\ref{tab:ratingScheme}.}
\label{tab:casesExamples}
\centering
\begin{tabular}{cc|c|c|c|c}
\multicolumn{2}{c|}{example} & \multicolumn{2}{c|}{\texttt{illustrative-1}} & \multicolumn{2}{c}{\texttt{illustrative-2}} \\
\multicolumn{2}{c|}{$\func$} & \multicolumn{2}{c|}{$\exp(x_1-2\,x_2^2+3\,x_3^3)$} & \multicolumn{2}{c}{$\frac{x_1}{x_1+0.2\,x_2^2}-2\,\frac{x_2}{x_2+0.3\,x_3^3}$}\\
\hline 
\hline 
 &$[\lb{x}_1,\ub{x}_1]$ & $[-0.3,0.2]$ & $[-0.198,0.177]$ &$[ 1.043,        1.535]$ & $[ 1.5 ,1.6]$\\
$\hyperRec$ &$[\lb{x}_2,\ub{x}_2]$ & $[-0.1,0.6]$ &  $[-0.473 ,         0.2]$&$ [0.6        ,1.969]$ & $[0.6,       1.1]$ \\
 & $[\lb{x}_3,\ub{x}_3]$ & $[-0.4,0.5]$ & $[ -0.392,0.39]$ & $[0.555 ,       0.772]$ & $[1.0,         1.6]$ \\
\hline
A &$[\lb{\lambda}_{\mathrm{A}},\ub{\lambda}_{\mathrm{A}}]$ &$[ -19.904,       37.004]$ & $[-15.767,19.27]$ &$[-43.934,27.391]$ &$[-45.014,       17.624]$ \\
G &$[\lb{\lambda}_{\mathrm{G}},\ub{\lambda}_{\mathrm{G}}]$ & $[-26.391,       38.587]$ &$[-15.767,18.443]$ &$[-44.907,27.391]$ &$[-40.725,       19.507]$\\
H &$[\lb{\lambda}_{\mathrm{H}},\ub{\lambda}_{\mathrm{H}}]$  &$[-20.597 ,      29.603]$ & $[-12.603,14.278]$ &$[-34.743,26.399]$ &$[-33.691,18.897]$\\
\hline
class & $\ub{\lambda}_{\mathrm{A}}$ & ($+$) & ($-$) &  ($\circ$) &  ($++$)\\
 & $\lb{\lambda}_{\mathrm{A}}$ & ($++$) & ($\circ$) & ($+$) & ($-$)
\end{tabular}
\end{table}

\subsection{Results}\label{subsec:Results}
Table \ref{tab:benchmark} summarizes the results obtained for the \numOfEx{} sample functions. The numbers listed in the columns labeled $\lb{\lambda}_\mathrm{A}$ state for how many of 100 randomly generated hyperrectangles the lower bounds calculated with the three methods fall into the classes $(-)$, $(\circ)$, $(+)$ and $(++)$ defined in Tab.~\ref{tab:ratingScheme}. The numbers listed in the columns labeled $\ub{\lambda}_\mathrm{A}$ state the corresponding results for the upper bounds. 

\begin{minipage}[t]{1.0\linewidth}
\small
\captionof{table}{Summary of results for the \numOfEx{} sample functions extracted from the COCONUT collection. 
For each example, lower and upper bounds on Hessian matrix eigenvalues were calculated with the three methods introduced in Lemmas \ref{lem:IntervalGershgorin}--\ref{lem:directArithmetic} for 100 random hyperrectangles.
Numbers state for how many out of the 100 random hyperrectangles the bounds belonged to the classes $(-)$, $(\circ)$, $(+)$ and $(++)$ defined in Tab. \ref{tab:ratingScheme}. 
Shaded cells highlight empty classes. Horizontal lines divide characteristic groups, e.g. examples with an empty class $(++)$ (rank $212\leq r \leq 854$).
The averages listed in the last row take all \numOfEx{} examples into account, including the ones not shown here.}
\label{tab:benchmark}  
\centering
\begin{tabular}{rcc||r|r|r|r||r|r|r|r} 
rank & \multicolumn{2}{c||}{example} & \multicolumn{4}{c||}{benchmark $\lb{\lambda}_{\mathrm{A}}$} & \multicolumn{4}{c}{benchmark $\ub{\lambda}_{\mathrm{A}}$} \\  
$r$ & name & $n$ & $(-)$ & $(\circ)$ & $(+)$ & $(++)$ & $(-)$ & $(\circ)$ & $(+)$ & $(++)$  \\ 
\hline 
\hline 
1  & \verb"box3-1"&  3& 2& 50& 32& 16& 0& 0& 2& 98\\
2  & \verb"box2-1"&  3& 0& 56& 26& 18& 0& 0& 6& 94\\
3  & \verb"cliff-1"&  2& 22& 54& 1& 23& 0& 10& 1& 89\\
4  & \verb"chaconn1-1"&  3& 29& 45& 0& 26& 0& 14& 0& 86\\
5  & \verb"chaconn2-1"&  3& 19& 63& 0& 18& 0& 11& 0& 89\\
6  & \verb"cb3-1"&  3& 18& 62& 0& 20& 0& 14& 0& 86\\
7  & \verb"polak6-1"&  5& 0& 0& 100& 0& 0& 0& 0& 100\\
8  & \verb"polak6-2"&  5& 0& 0& 100& 0& 0& 0& 0& 100\\
9  & \verb"polak6-3"&  5& 0& 0& 100& 0& 0& 0& 0& 100\\
10  & \verb"polak6-4"&  5& 0& 0& 100& 0& 0& 0& 0& 100\\
11  & \verb"growth-1"&  3& 0& 0& 96& 4& 0& 0& 4& 96\\
12  & \verb"alsotame-1"&  2& 0& 0& 100& 0& 0& 0& 0& 100\\
13  & \verb"vardim-1"&  10& 0& 0& 0& 100& 0& 0& 100& 0\\
14  & \verb"vardim-2"&  10& 0& 0& 100& 0& 0& 0& 0& 100\\
15  & \verb"alsotame-2"&  2& 0& 0& 100& 0& 0& 0& 0& 100\\
16  & \verb"brownden-1"&  4& 1& 0& 99& 0& 0& 0& 0& 100\\
17  & \verb"price-1"&  2& 3& 0& 97& 0& 0& 0& 0& 100\\
18  & \verb"vanderm1-10"&  10& 75& 0& 25& 0& 0& 0& 0& 100\\
19  & \verb"ex8_1_7-1"&  5& 99& 0& 1& 0& 0& 0& 0& 100\\
20  & \verb"hs026-2"&  3& 0& 100& 0& 0& 0& 0& 0& 100\\
\vdots \\
209  & \verb"ex14_1_7-5"&  10& 4& 0& 95& 1& 25& 0& 75& 0\\
210  & \verb"ex14_1_7-9"&  10& 25& 0& 75& 0& 4& 0& 95& 1\\
211  & \verb"nonmsqrt-1"&  9& 100& 0& 0& 0& 96& 0& 3& 1\\
\hline
212  & \verb"brkmcc-1"&  2& 0& 0& 100& \cellcolor[gray]{0.8} 0& 0& 0& 100& \cellcolor[gray]{0.8} 0\\
213  & \verb"ship-15"&  10& 1& 0& 99& \cellcolor[gray]{0.8} 0& 1& 0& 99& \cellcolor[gray]{0.8} 0\\
\vdots & & \multicolumn{1}{c}{}& \multicolumn{1}{r}{}  & \multicolumn{1}{r}{}  &  \multicolumn{1}{r}{}  & \multicolumn{1}{r}{\cellcolor[gray]{0.8} \vdots}  & \multicolumn{1}{r}{}  & \multicolumn{1}{r}{}  & \multicolumn{1}{r}{}   & \multicolumn{1}{r}{\cellcolor[gray]{0.8} \vdots}  \\
852  & \verb"butcher-4"&  7& 100& 0& 0& \cellcolor[gray]{0.8} 0& 99& 0& 1& \cellcolor[gray]{0.8} 0\\
853  & \verb"i5-3"&  10& 100& 0& 0& \cellcolor[gray]{0.8} 0& 99& 0& 1& \cellcolor[gray]{0.8} 0\\
854  & \verb"cohn3-1"&  4& 100& 0& 0& \cellcolor[gray]{0.8} 0& 99& 0& 1& \cellcolor[gray]{0.8} 0\\
\hline
855  & \verb"ex4_1_8-1"&  2& 0& 100& \cellcolor[gray]{0.8} 0& \cellcolor[gray]{0.8} 0& 0& 100& \cellcolor[gray]{0.8} 0& \cellcolor[gray]{0.8} 0\\
856  & \verb"sample-3"&  4& 0& 100& \cellcolor[gray]{0.8} 0& \cellcolor[gray]{0.8} 0& 0& 100& \cellcolor[gray]{0.8} 0& \cellcolor[gray]{0.8} 0\\
\vdots & & \multicolumn{1}{c}{}& \multicolumn{1}{r}{}  & \multicolumn{1}{r}{}  & \multicolumn{1}{r}{\cellcolor[gray]{0.8} \vdots}  & \multicolumn{1}{r}{\cellcolor[gray]{0.8} \vdots}  & \multicolumn{1}{r}{}  & \multicolumn{1}{r}{}  & \multicolumn{1}{r}{\cellcolor[gray]{0.8} \vdots}  & \multicolumn{1}{r}{\cellcolor[gray]{0.8} \vdots}   \\
1391  & \verb"womflet-1"&  3& 100& 0& \cellcolor[gray]{0.8} 0& \cellcolor[gray]{0.8} 0& 97& 3& \cellcolor[gray]{0.8} 0& \cellcolor[gray]{0.8} 0\\
1392  & \verb"reimer5-2"&  5& 98& 2& \cellcolor[gray]{0.8} 0& \cellcolor[gray]{0.8} 0& 99& 1& \cellcolor[gray]{0.8} 0& \cellcolor[gray]{0.8} 0\\
1393  & \verb"reimer5-5"&  5& 98& 2& \cellcolor[gray]{0.8} 0& \cellcolor[gray]{0.8} 0& 99& 1& \cellcolor[gray]{0.8} 0& \cellcolor[gray]{0.8} 0\\
\hline
1394  & \verb"ex7_2_9-4"&  10& 100& \cellcolor[gray]{0.8} 0& \cellcolor[gray]{0.8} 0& \cellcolor[gray]{0.8} 0& 100& \cellcolor[gray]{0.8} 0& \cellcolor[gray]{0.8} 0& \cellcolor[gray]{0.8} 0\\
1395  & \verb"ex7_2_9-2"&  10& 100& \cellcolor[gray]{0.8} 0& \cellcolor[gray]{0.8} 0& \cellcolor[gray]{0.8} 0& 100& \cellcolor[gray]{0.8} 0& \cellcolor[gray]{0.8} 0& \cellcolor[gray]{0.8} 0\\
\vdots & & \multicolumn{1}{c}{}& \multicolumn{1}{r}{\vdots}  & \multicolumn{1}{r}{\cellcolor[gray]{0.8} \vdots}  & \multicolumn{1}{r}{\cellcolor[gray]{0.8} \vdots}  & \multicolumn{1}{r}{\cellcolor[gray]{0.8} \vdots}  & \multicolumn{1}{r}{\vdots}  & \multicolumn{1}{r}{\cellcolor[gray]{0.8} \vdots}  & \multicolumn{1}{r}{\cellcolor[gray]{0.8} \vdots}  & \multicolumn{1}{r}{\cellcolor[gray]{0.8} \vdots}   \\
1518  & \verb"cohn2-2"&  4& 100& \cellcolor[gray]{0.8} 0& \cellcolor[gray]{0.8} 0& \cellcolor[gray]{0.8} 0& 100& \cellcolor[gray]{0.8} 0& \cellcolor[gray]{0.8} 0& \cellcolor[gray]{0.8} 0\\
1519  & \verb"cohn2-3"&  4& 100& \cellcolor[gray]{0.8} 0& \cellcolor[gray]{0.8} 0& \cellcolor[gray]{0.8} 0& 100& \cellcolor[gray]{0.8} 0& \cellcolor[gray]{0.8} 0& \cellcolor[gray]{0.8} 0\\
1520  & \verb"cohn2-4"&  4& 100& \cellcolor[gray]{0.8} 0& \cellcolor[gray]{0.8} 0& \cellcolor[gray]{0.8} 0& 100& \cellcolor[gray]{0.8} 0& \cellcolor[gray]{0.8} 0& \cellcolor[gray]{0.8} 0\\
1521  & \verb"boon-2"&  6& 100& \cellcolor[gray]{0.8} 0& \cellcolor[gray]{0.8} 0& \cellcolor[gray]{0.8} 0& 100& \cellcolor[gray]{0.8} 0& \cellcolor[gray]{0.8} 0& \cellcolor[gray]{0.8} 0\\
1522  & \verb"boon-4"&  6& 100& \cellcolor[gray]{0.8} 0& \cellcolor[gray]{0.8} 0& \cellcolor[gray]{0.8} 0& 100& \cellcolor[gray]{0.8} 0& \cellcolor[gray]{0.8} 0& \cellcolor[gray]{0.8} 0\\
\hline 
\multicolumn{3}{c||}{arithmetic average}& 24.17& 60.85& 12.48& 2.50& 23.52& 60.33& 10.63& 5.52\\
\end{tabular}
\end{minipage}

The examples are ranked in Tab.~\ref{tab:benchmark} by, loosely speaking, the quality of the bounds found with the eigenvalue arithmetic. More precisely, an example is ranked the higher, the higher the sum of the figures in its two $(++)$ columns. If this sum is equal for several examples, they are sorted according to the sum of the figures in their two $(+)$ columns.  Subsequently, the sums of the two $(\circ)$ columns and the two $(-)$ columns are used for the ranking whenever necessary.

Table \ref{tab:benchmark} shows the 20 best and 5 worst rated examples and characteristic ranks in between. Ranks \rankFirstBreakUpper{}$-$\rankFirstBreakLower{}, for example, are shown, because they mark the boundary between those $\func$ for which some eigenvalue bounds still fall into class $(++)$, and the highest ranking examples for which class $(++)$ no longer occurs. These transitions in the ranking are marked with horizontal lines and shaded areas. 

Just as for the illustrative examples given in Tab.~\ref{tab:casesExamples}, all classes $(-)$, $(\circ)$, $(+)$, $(++)$ occur for the sample functions extracted from the COCONUT collection. 
In particular there exist cases for which the eigenvalue arithmetic provides tighter bounds than the tight bounds for the interval Hessian. An analysis of the data given in Tab.~\ref{tab:benchmark} reveals that $\relNumPlusPlusLower\%$ and $\relNumPlusPlusUpper\%$ of the examples belong to the $(++)$ class for the lower bound and upper eigenvalue bound, respectively. 
Furthermore, in $\relNumPlusLower\%$ of the cases the eigenvalue arithmetic provides tighter lower bounds than Gershgorin's circle criterion applied to the interval Hessian matrix. In $\relNumPlusUpper\%$ of the cases the upper bound from the eigenvalue arithmetic is tighter than the upper Gershgorin bound.   
In $\relNumCircLower\%$ ($\relNumCircUpper\%$) of the cases the lower (upper) bounds from the eigenvalue arithmetic and those from Gershgorin's circle criterion are equal. We stress, however, that the Gershgorin bounds outperform those of the eigenvalue arithmetic in $\relNumMinusLower\%$ (lower bound) and $\relNumMinusUpper\%$ (upper bound) of the cases, respectively. 
Overall, there exist \numAllBetter{} examples for which the eigenvalue arithmetic provides better bounds than Gershgorin's circle criterion for all random boxes (e.g. $7\leq r\leq 15$ or $r=212$). On the other hand, there exist  \numAllWorse{} examples for which the arithmetic results in less tight spectral bounds for all random boxes ($r\geq\boundAllWorse$).

For \numAtLeastOneBetter{} out of the \numOfEx{} examples (i.e., $\relNumAtLeastOneBetter\%$) the eigenvalue arithmetic provides tighter bounds than Gershgorin for at least one of the random boxes.
Finally, we note that the ranking $r$ does not correlate with the dimension $n$, i.e., we find both low and high values of $n$ in any part of the ranking shown in Tab.~\ref{tab:benchmark}. The dependency on $n$ is analyzed in more detail at the end of this section with Tab.~\ref{tab:dimension}.

\begin{table}[hp] 
\caption{Computational effort for the evaluation of eigenvalue bounds for the examples taken from the COCONUT-benchmark.
$\NA$ and $\NGersh$ denote the total number of operations necessary to calculate
$[\lambda_{\mathrm{A}}]$ and $[\lambda_{\mathrm{G}}]$, respectively. 
$\DeltaNA$ denotes the additional number of operations needed to calculate $[\lambda_{\mathrm{A}}]$ assuming the Gershgorin bounds have already been computed. 
The examples are listed in the same order as in Tab.~\ref{tab:benchmark}. 
Shaded cells highlight minimal and maximal values for the ratios $\frac{\NA}{\NGersh}$ and $\frac{\DeltaNA}{\NGersh}$, respectively.  
The averages stated in the last row take all \numOfEx{} examples into account, including those not shown here.} 
\label{tab:complexity} 
\begin{tabular}{rcc||rrr|rr} 
\multicolumn{3}{c||}{example} & \multicolumn{3}{c|}{abs. complexity} & \multicolumn{2}{c}{rel. complexity (\%)} \\  
$r$ & name & $n$ & $\NA$ & $\NGersh$  & $\DeltaNA$ & $\frac{\NA}{\NGersh}$ & $\frac{\DeltaNA}{\NGersh}$ \\ 
\hline 
\hline 
1  & \verb"box3-1" & 3& 2270& 5798& 808& 39.15& 13.94\\
2  & \verb"box2-1" & 3& 2292& 5820& 808& 39.38& 13.88\\
3  & \verb"cliff-1" & 2& 159& 289& 56& 55.02& 19.38\\
4  & \verb"chaconn1-1" & 3& 87& 242& 29& 35.95& 11.98\\
5  & \verb"chaconn2-1" & 3& 87& 242& 29& 35.95& 11.98\\
6  & \verb"cb3-1" & 3& 87& 242& 29& 35.95& 11.98\\
7  & \verb"polak6-1" & 5& 1302& 5526& 456& 23.56& 8.25\\
8  & \verb"polak6-2" & 5& 1300& 5524& 456& 23.53& 8.25\\
9  & \verb"polak6-3" & 5& 1302& 5526& 456& 23.56& 8.25\\
10  & \verb"polak6-4" & 5& 1302& 5526& 456& 23.56& 8.25\\
11  & \verb"growth-1" & 3& 3470& 8132& 1546& 42.67& 19.01\\
12  & \verb"alsotame-1" & 2& 59& 120& 23& 49.17& 19.17\\
\hline
13  & \verb"vardim-1" & 10& 2828& 21232& 721& \cellcolor[gray]{0.8}13.32& 3.40\\
\hline
14  & \verb"vardim-2" & 10& 3118& 22818& 745& 13.66& \cellcolor[gray]{0.8}3.26\\
\hline
15  & \verb"alsotame-2" & 2& 61& 122& 23& 50.00& 18.85\\
16  & \verb"brownden-1" & 4& 5696& 18610& 1896& 30.61& 10.19\\
17  & \verb"price-1" & 2& 423& 769& 192& 55.01& 24.97\\
18  & \verb"vanderm1-10" & 10& 18188& 174048& 6618& 10.45& 3.80\\
19  & \verb"ex8_1_7-1" & 5& 685& 2920& 229& 23.46& 7.84\\
20  & \verb"hs026-2" & 3& 198& 530& 75& 37.36& 14.15\\
\vdots \\
\hline
 161  & \verb"gold-1" & 2& 961& 1600& 447& \cellcolor[gray]{0.8}60.06& 27.94\\
\hline
\vdots \\
\hline
838  & \verb"sendra-1" & 2& 599& 1038& 311& 57.71& \cellcolor[gray]{0.8}29.96\\
\hline
\vdots \\
1518  & \verb"cohn2-2" & 4& 3132& 9296& 1514& 33.69& 16.29\\
1519  & \verb"cohn2-3" & 4& 3132& 9296& 1514& 33.69& 16.29\\
1520  & \verb"cohn2-4" & 4& 5067& 14578& 2463& 34.76& 16.90\\
1521  & \verb"boon-2" & 6& 718& 3298& 350& 21.77& 10.61\\
1522  & \verb"boon-4" & 6& 718& 3298& 350& 21.77& 10.61\\
\hline 
\multicolumn{3}{c||}{arithmetic average} & \multicolumn{3}{c|}{\it irrelevant}& 37.40& 17.35\\
\end{tabular} 
\end{table} 

We discussed in Sect.~\ref{subSec:Complexity} that the compared methods do not only differ with respect to the tightness of eigenvalue bounds, but also with respect to their computational complexity. 
Table \ref{tab:complexity} lists the operation numbers for the examples from Tab.~\ref{tab:benchmark}.
All figures in Tab.~\ref{tab:complexity} are based on the total number of operations needed for the respective method. 
Specifically, $\NA$ denotes the total number of operations for calculating $[\lambda_\mathrm{A}]$ 
with the eigenvalue arithmetic, 
including the operations for the intermediate results $[y_k]$ and $[y_k^\prime]$. 
$\NGersh$ denotes the total number of operations for calculating $[\lambda_\mathrm{G}]$ 
with Gershgorin's circle criterion, 
including the operations for the intermediate results $[y_k]$, $[y_k^\prime]$ and $[y_k^{\prime\prime}]$.  
We also list $\DeltaNA$ defined in (\ref{eq:additionalComplexity}), i.e.\ the additional effort to calculate $[\lambda_\mathrm{A}]$, if $[\lambda_\mathrm{G}]$ and its intermediate results $[y_k]$ and $[y_k^\prime]$ have been determined. 
As predicted by relation (\ref{eq:complexityRelationAG}), the eigenvalue arithmetic always requires fewer operations than the interval variant of Gershgorin's circle criterion. 
On average the computational effort for the eigenvalue arithmetic amounts to $\meanComplexityNAOverNG\%$ of that of applying Gershgorin's circle criterion to the interval Hessian, where this figure ranges from $\minComplexityNAOverNG\%$ (example $r=\rankMinComplexityNAOverNG$ in Tab.~\ref{tab:complexity}) to $\maxComplexityNAOverNG\%$ (example $r=\rankMaxComplexityNAOverNG$). The average additional effort $\DeltaNA$ for the eigenvalue arithmetic amounts to $\meanComplexityDeltaNAOverNG\%$, with a minimum and maximum of $\minComplexityDeltaNAOverNG\%$ (example $r=\rankMinComplexityDeltaNAOverNG$) and $\maxComplexityDeltaNAOverNG\%$ (example $r=\rankMaxComplexityDeltaNAOverNG$).

We stress that $\NA$, $\NGersh$, $\DeltaNA$ do not depend on the particular hyperrectangle $\hyperRec$, but can be determined for any function $\func$ before eigenvalue bounds are actually calculated. It may therefore be an option to determine these operation counts beforehand and to decide which method to use. This may be an option in applications in which eigenvalue bounds need to be calculated for the same function $\func$ for many $\hyperRec$ such as branch-and-bound global optimization.

While we did not recognize a dependency of the tightness of the bounds on $n$, the ratios  $\NA / \NGersh$ and $\DeltaNA/\NGersh$ clearly depend on $n$. This was anticipated in the discussion of the complexity classes in Sect.~\ref{subSec:Complexity}. Table~\ref{tab:dimension} shows that the relative number of operations for the eigenvalue arithmetic $\NA / \NGersh$ improves from about $\approxMeanComplexityNAOverNGForNTwo\%$ for $n=2$ to about $\approxMeanComplexityNAOverNGForNTen\%$ for $n= 10$.
Similarly, 
the additional effort for the eigenvalue arithmetic decreases from about $\approxMeanComplexityDeltaNAOverNGForNTwo\%$ for $n= 2$ to about $\approxMeanComplexityDeltaNAOverNGForNTen\%$ for $n= 10$. Note that examples with $n=1$ would yield $\NA / \NGersh=1=100\%$.

Finally, we note that the large number of examples $\func$ for $n=3$ results from the COCONUT optimization problem \texttt{oet2}. We did not omit any of these $\func$ in order not to introduce bias. 

\begin{table}[h] 
\caption{Aggregation of results by dimension $n$.} 
\label{tab:dimension} 
\centering 
\begin{tabular}{rr|rrrr|rr|rr} 
& \multicolumn{1}{c|}{num. of} &  \multicolumn{4}{c|}{mean: (\%)} & \multicolumn{2}{c|}{mean: (\%)} & \multicolumn{2}{c}{std: (\%)}\\ 
$n$ & examples & \multicolumn{1}{c}{$(-)$} & \multicolumn{1}{c}{$(\circ)$} & \multicolumn{1}{c}{$(+)$} & \multicolumn{1}{c|}{$(++)$} & $\frac{\NA}{\NGersh}$ & $\frac{\DeltaNA}{\NGersh}$ & $\frac{\NA}{\NGersh}$ & $\frac{\DeltaNA}{\NGersh}$\\ 
\hline 
\hline 
2 & 62& 57.89& 15.47& 14.68& 11.97& \cellcolor[gray]{0.1}\color{white}54.16& \cellcolor[gray]{0.1}\color{white}24.41& 3.01& 3.17\\
3 & 1078& 10.88& 79.32& 8.79& 1.01& \cellcolor[gray]{0.2}\color{white}41.70& \cellcolor[gray]{0.2}\color{white}19.81& 1.04& 1.21\\
4 & 67& 61.29& 19.34& 8.45& 10.92& \cellcolor[gray]{0.3}\color{white}31.90& \cellcolor[gray]{0.3}\color{white}13.75& 2.76& 2.42\\
5 & 88& 56.86& 15.85& 12.48& 14.81& \cellcolor[gray]{0.4}\color{white}25.34& \cellcolor[gray]{0.4}\color{white}10.33& 3.14& 2.13\\
6 & 95& 35.05& 14.21& 36.88& 13.86& \cellcolor[gray]{0.5}23.08& \cellcolor[gray]{0.5}8.97& 2.12& 1.47\\
7 & 27& 65.81& 34.17& 0.02& 0.00& \cellcolor[gray]{0.6}18.98& \cellcolor[gray]{0.6}7.44& 1.76& 1.60\\
8 & 15& 94.23& 4.50& 1.27& 0.00& \cellcolor[gray]{0.7}17.85& \cellcolor[gray]{0.7}7.83& 0.99& 0.82\\
9 & 24& 65.60& 4.21& 18.71& 11.48& \cellcolor[gray]{0.8}14.48& \cellcolor[gray]{0.8}6.15& 2.60& 1.80\\
10 & 66& 57.06& 9.34& 23.74& 9.86& \cellcolor[gray]{0.9}14.11& \cellcolor[gray]{0.9}5.90& 2.22& 1.72\\
\hline 
all & 1522& 23.84& 60.59& 11.55& 4.01& 37.40& 17.35& 9.56& 5.17
\end{tabular} 
\end{table}

\section{Conclusion and Outlook}
\label{sec:Conclusion}
Our numerical experiments corroborate that 
the eigenvalue arithmetic always requires fewer operations than Gershgorin's circle criterion. While this result has been established by comparing the complexity classes of the two methods (see Sect.\ \ref{subSec:Complexity} and \cite{Monnigmann2011}), it was analyzed quantitatively with a large set of examples here for the first time. Specifically, 
$\minComplexityNAOverNG\%$ ($n= 10$)  to $\maxComplexityNAOverNG\%$ ($n= 2$)
of the number of operations of the Gershgorin based approach are necessary for the eigenvalue method. The average over all examples for all $n$ amounts to $\meanComplexityNAOverNG\%$. As anticipated in the complexity analysis in Sect.~\ref{subSec:Complexity}, the eigenvalue method benefits from its favorable complexity as $n$ increases (see Tab.~\ref{tab:complexity} for details). 
We recall that a comparison to the computational effort of Hertz and Rohn's method is not reasonable, since Hertz and Rohn's method belongs to a very different complexity class (see Sect.\ \ref{subSec:Complexity}).  

Gershgorin's circle criterion provides tighter lower (upper) bounds in $\relNumMinusLower\%$ ($\relNumMinusUpper\%$) of the examples. In $\relNumCircLower\%$ ($\relNumCircUpper\%$) of the cases the lower (upper) bounds from both methods are equal. 
In $\relNumPlusLower\%$ ($\relNumPlusUpper\%$) of the examples the eigenvalue arithmetic results in tighter lower (upper) bounds than the Gershgorin based approach. Finally, our tests reveal that the number of cases in which the eigenvalue arithmetic results in tighter bounds than the tight bounds of the interval Hessian, which are obtained with Hertz and Rohn's method, is small ($\relNumPlusPlusLower\%$ and $\relNumPlusPlusUpper\%$ for lower and upper bounds, respectively). On the other hand, these figures indicate that these cases are not anecdotal or constructed, but they appear in global optimization problems. 

The eigenvalue arithmetic provides a tighter lower or upper bound than Gershgorin's circle criterion for at least one random box in \relNumAtLeastOneBetter\% of the examples, where these occurencies are not correlated with $n$. This figure suggests to combine the two methods. We claim the eigenvalue arithmetic can be applied at an attractive additional cost for, say, $n> 5$, whenever the Gershgorin bounds have already been calculated, since both methods involve the same intermediate quantities ($[y_k]$ and $[y^\prime_k]$, see Sect.~\ref{sec:Arithmetic}). Specifically, the additional effort for applying the eigenvalue method after the intermediate quantities have been calculated in the Gershgorin based approach ranges from $\approxMeanComplexityDeltaNAOverNGForNTwo\%$ for $n= 2$ to about $\approxMeanComplexityDeltaNAOverNGForNTen\%$ for $n= 10$ (see Tab.~\ref{tab:dimension}). This figure decreases for increasing $n$ as anticipated from the abstract complexity analysis in Sect.~\ref{subSec:Complexity}. Note that this combination of Gershgorin's circle criterion and the eigenvalue method will provide tighter bounds than Hertz and Rohn's method for the interval Hessian whenever the eigenvalue method does.

\end{document}